\documentclass[12pt]{article}
\usepackage{amsfonts}

\usepackage{latexsym}
\usepackage{amssymb}
\usepackage{amsmath}
\usepackage{amsthm}
\usepackage[mathscr]{eucal}
\usepackage{graphicx}
\usepackage{hyperref}
\usepackage{caption}
\usepackage{subcaption}

\newtheorem{Lemma}{Lemma}

\newtheorem{Theorem}[Lemma]{Theorem}

\renewcommand{\qed}{\hfill{\ \ \rule{2mm}{2mm}} \vspace{0.2in}}

\newcommand{\ind}{1\hspace{-2.3mm}{1}}

\begin{document}

\title{Traveling salesman problem across dense cities}
\author{ \textbf{Ghurumuruhan Ganesan}
\thanks{E-Mail: \texttt{gganesan82@gmail.com} } \\
\ \\
New York University, Abu Dhabi }
\date{}
\maketitle

\begin{abstract}
Consider~\(n\) nodes~\(\{X_i\}_{1 \leq i \leq n}\) distributed independently across~\(N\) cities
contained with the unit square~\(S\) according to a distribution~\(f.\)
Each city is modelled as an~\(r_n \times r_n\) square contained within~\(S\)
and let~\(TSPC_n\) denote the length of the minimum length cycle containing
all the~\(n\) nodes, corresponding to the traveling salesman
problem (TSP). We obtain variance estimates for~\(TSPC_n\)
and prove that if the cities are well-connected and densely populated
in a certain sense, then~\(TSPC_n\) appropriately centred and scaled
converges to zero in probability. We also obtain
large deviation type estimates for~\(TSPC_n.\) Using the proof techniques,
we alternately obtain corresponding results for the length~\(TSP_n\)
of the minimum length cycle in the unconstrained case, when the nodes
are independently distributed throughout
the unit square~\(S.\)



\vspace{0.1in} \noindent \textbf{Key words:} Traveling salesman problem, dense cities.

\vspace{0.1in} \noindent \textbf{AMS 2000 Subject Classification:} Primary:
60J10, 60K35; Secondary: 60C05, 62E10, 90B15, 91D30.
\end{abstract}

\bigskip

\setcounter{equation}{0}
\renewcommand\theequation{\thesection.\arabic{equation}}
\section{Introduction} \label{intro}
The Traveling Salesman Problem (TSP) is the study of finding the minimum weight cycle containing all
the nodes of a graph where each edge is assigned a certain weight.
In this paper, we consider the case of random Euclidean TSP, henceforth referred to simply
as TSP, where the nodes are distributed randomly across the unit square~\(S\) with origin as centre.
The weight of an edge between two nodes is the Euclidean distance between them
and the goal is to find the cycle of shortest length containing all the nodes.
For more material on the TSP, we refer to the books by Gutin and Punnen (2006),
Cook (2011) and references therein.

The analytical study of the random TSP problem originated in Beardwood et al (1959).
The main result there is that if~\(n\) nodes are randomly and uniformly distributed
across the unit square~\(S,\)
then with high probability (i.e., with probability converging to one as~\(n \rightarrow \infty\)),
the length~\(TSP_n\) of the minimum length spanning cycle grows roughly as~\(\beta \sqrt{n}\)
for some constant~\(\beta > 0.\) Equivalently,~\(TSP_n\) appropriately
scaled and centred converges to zero a.s.\ and in mean as~\(n \rightarrow \infty.\)
Subadditive ergodic type theorems are used for obtaining the convergence results
and for a comprehensive survey, we refer to Steele~(1981, 1993).

Since then there has been a lot of work focused on obtaining better bounds
for the constant~\(\beta > 0.\) Beardwood et al originally established
that~\(0.625 \leq \beta \leq 0.922.\) Recently, Steinerberger (2015)
has obtained slightly improved bounds by estimating the probability
of certain configurations that are avoided by the optimal cycle.

Because of its practical importance, there has also been a lot of work
devoted to obtaining optimal and near optimal algorithms for obtaining the
minimum length cycle. Arora (1998), Vazirani (2001), Karpinski et al (2015)
develop and analyse polynomial time approximation
schemes (PTAS) that determine near minimal spanning cycles for large vertex sets.
Snyder and Daskin (2006) have used genetic algorithms to provide heuristic solutions
for the generalized TSP problem, where the nodes are split into clusters and the objective
is to find a minimum cost tour passing through exactly one node from each cluster.
Recently, Pintea et al (2017) have proposed solutions to the generalized TSP problem
using Ant algorithms.

The analytical literature above mainly consider nodes distributed in regular shapes like unit squares or circles.
In this paper, we consider a slightly different scenario where cities (modelled as small squares) are spread
across the unit square each containing a subset of the nodes. The cities are not necessarily regularly spaced and therefore
the usual subadditive techniques to determine the convergence of TSP are not directly applicable here.
Instead, we use approximation methods to find sharp upper and lower bounds
for the optimal minimum spanning cycle and indirectly deduce convergence properties as the size of
the vertex set~\(n \rightarrow \infty.\)

\subsection*{Model Description}
\subsubsection*{Structure of the cities}
\begin{figure}[tbp]
\centering
\includegraphics[width=2in, trim= 100 180 50 100, clip=true]{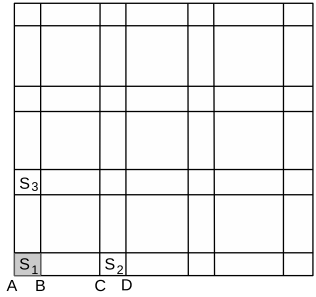}
\caption{Tiling~\(S\) into~\(r_n \times r_n\) squares with an inter-square distance of~\(s_n.\)}
\label{sq_plc}
\end{figure}

For integer~\(n \geq 1,\) let~\(r_n\) and~\(s_n\) be real numbers such that~\(\frac{1-r_n}{r_n+s_n}\)
is an integer. Tile the unit square~\(S\) regularly into~\(r_n \times r_n\) size squares
in such a way that the distance between any two squares is at least~\(s_n\)
as shown in Figure~\ref{sq_plc}. In Figure~\ref{sq_plc}, the grey square is of
size~\(r_n  \times r_n,\)  the segment~\(AB\) has length~\(r_n\) and the segment~\(BC\)
has length~\(s_n.\) The~\(r_n \times r_n\) squares are
called \emph{cities} and the term~\(s_n\) denotes the
\emph{intercity distance}.

Label the~\(r_n \times r_n\) squares (cities) as~\(\{S_l\}\)
and identifying the centres of the squares~\(\{S_l\}\) with vertices in~\(\mathbb{Z}^2,\)
we obtain a corresponding subset of vertices~\(\{z_l\} \subset \mathbb{Z}^2.\)
For example, in Figure~\ref{sq_plc}, identify the centre of the square labelled~\(S_1\)
with~\((0,0),\) the centre of~\(S_2\) with~\((1,0),\) the centre of~\(S_3\) with~\((0,1)\)
and so on. Two vertices~\(z_1 = (x_1,y_1)\) and~\(z_2 = (x_2,y_2)\) are \emph{adjacent}
and connected by an edge if~\(|x_1-x_2| + |y_1 - y_2| = 1.\) 

Fix~\(N = N(n)\) cities~\(\{S_{j_1},\ldots,S_{j_N}\}\) and let~\(\{z_{j_1},\ldots,z_{j_N}\}\)
be the vertices in~\(\mathbb{Z}^2\) corresponding to
the centres of~\(\{S_{j_i}\}.\) We say that the cities~\(\{S_{j_1},\ldots,S_{j_N}\}\)
are \emph{well-connected} if the corresponding set of vertices~\(\{z_{j_i}\}\)
form a connected subgraph of~\(\mathbb{Z}^2.\) Henceforth, we
assume that~\(\{S_{j_1},\ldots,S_{j_N}\}\) are well-connected
and without loss of generality denote~\(S_{j_i}\) by~\(S_i\)
for~\(1 \leq i \leq N.\)

\subsubsection*{Nodes in the cities}
Let~\(f\) be any density on the unit square~\(S\) satisfying the following conditions:\\
There are constants~\(0 < \epsilon_1 \leq \epsilon_2 < \infty\) such that
\begin{equation}\label{f_eq}
\epsilon_1 \leq \inf_{x \in S} f(x) \leq \sup_{x \in S} f(x) \leq \epsilon_2
\end{equation}
and
\begin{equation}\label{f_eq_tot}
\int_{x \in S} f(x) dx = 1.
\end{equation}
Define the density~\(g_N(.)\) on the~\(N\) cities~\(\bigcup_{1 \leq i \leq N} S_i\)
as
\begin{equation}\label{gn_def}
g_N(x) = \frac{f(x)}{\int_{\cup_{1 \leq j \leq N}S_j} f_j(x) dx}
\end{equation}
for all~\(x \in \bigcup_{1 \leq j \leq N} S_j.\)


Let~\(X_1,X_2,\ldots,X_n\) be~\(n\) nodes independently
and identically distributed (i.i.d.)
in the~\(N\) cities~\(\{S_j\}_{1 \leq j \leq N},\) each according to
the density \(g_N.\) Define the vector~\((X_1,\ldots,X_n)\) on the probability
space \((\Omega_X, {\cal F}_X, \mathbb{P}).\) Let~\(K_n = K(X_1,\ldots,X_n)\) be the complete graph whose edges are
obtained by connecting each pair of nodes~\(X_i\) and~\(X_j\) by the straight line segment~\((X_i,X_j)\)
with~\(X_i\) and~\(X_j\) as endvertices. The line segment~\((X_i,X_j)\)
is the edge between the nodes~\(X_i\) and\(X_j\)
and~\(d(X_i,X_j)\) denotes the (Euclidean) length of the edge~\((X_i,X_j).\)


A cycle~\({\cal C}  = (Y_1,Y_2,\ldots,Y_t,Y_1)\) is a subgraph of~\(K_n\) with vertex set
\(\{Y_{j}\}_{1 \leq j \leq t} \subset \{X_i\}\)
and edge set~\(\{(Y_{j},Y_{{j+1}})\}_{1 \leq j \leq t-1} \cup (Y_t,Y_1).\)
The length of~\({\cal C}\) is defined as the sum of the lengths of the edges in~\({\cal C};\) i.e.,
\begin{equation}\label{len_cyc_def}
L({\cal C}) = \sum_{i=1}^{t-1} d(Y_i,Y_{i+1}) + d(Y_{t},Y_1) = \frac{1}{2} \sum_{i=1}^{t} l(Y_i,{\cal C}),
\end{equation}
where~\(l(Y_1,{\cal C}) = d(Y_1,Y_2) + d(Y_1,Y_t), l(Y_t,{\cal C}) = d(Y_{t},Y_1) + d(Y_t,Y_{t-1})\) and for~\(2 \leq i \leq t,\)
\[l(Y_i,{\cal C}) = d(Y_i,Y_{i-1}) + d(Y_i,Y_{i+1})\] is the sum of the length of the (two) edges in~\({\cal C}\) containing~\(Y_i\)
as an endvertex.
The cycle~\({\cal C}\) is said to be a \emph{spanning cycle} if~\({\cal C}\)
contains all the nodes~\(\{X_k\}_{1 \leq k \leq n}.\) Let~\({\cal C}_n\)
be a spanning cycle satisfying
\begin{equation}\label{min_weight_cycle}
TSPC_n = L({\cal C}_n) := \min_{{\cal C}} L({\cal C}),
\end{equation}
where the minimum is taken over all spanning cycles~\({\cal C}.\)
If there is more than one choice for~\({\cal C}_n,\)
choose one according to a deterministic rule. The cycle~\({\cal C}_n\)
is defined to the \emph{minimum spanning cycle} with corresponding length~\(TSPC_n.\)

Letting
\begin{equation}\label{bn_def}
b_n := r_n \sqrt{nN},
\end{equation}
we have the following result.
\begin{Theorem}\label{tsp_thm}
Suppose~\(r_n,s_n\) and~\(N = N(n)\) satisfy
\begin{equation}\label{N_est}
r_n^2 \geq \frac{M\log{n}}{n}, \frac{n}{N^2} \longrightarrow 0 \text{ and }\frac{N s_n}{b_n} \longrightarrow 0
\end{equation}
as~\(n \rightarrow \infty,\) for some constant~\(M >0.\) If~\(M =M(\epsilon_1,\epsilon_2) > 0\) is large, then
\begin{equation} \label{conv_tsp_prob}
\frac{1}{b_n}\left(TSPC_n - \mathbb{E}TSPC_n\right) \longrightarrow 0 \text{ in probability}
\end{equation}
as~\(n \rightarrow \infty.\) In addition, there are positive constants~\(\{\theta_i\}_{1 \leq i \leq 6}\) such that
\begin{equation}\label{exp_tspc}
\theta_1 b_n \leq \mathbb{E}TSPC_n \leq \theta_2 b_n,
\end{equation}
\begin{equation}\label{eq_tsp2}
\mathbb{P}\left(TSPC_n \geq  \theta_3 b_n \right) \geq 1 - e^{-\theta_4 N}
\end{equation}
and
\begin{equation}\label{eq_tsp3}
\mathbb{P}\left(TSPC_n \leq  \theta_5 b_n \right) \geq 1 - \exp\left(-\theta_6 \frac{n}{N}\right)
\end{equation}
for all~\(n\) large. 
\end{Theorem}
In words, if the cities are wide and dense enough, then the centred and scaled minimum length
of the traveling salesman cycle converges to zero in probability. 

\subsection*{Unconstrained TSP}
There are~\(n\) nodes~\(\{X_i\}_{1 \leq i \leq n}\)
independently distributed in the unit square~\(S\) each according
to the distribution~\(f\) satisfying~(\ref{f_eq}).
As in~(\ref{min_weight_cycle}), let~\(TSP_n\) be the length
of the minimum spanning cycle containing all the nodes~\(\{X_i\}_{1 \leq i \leq n}.\)

Beardwood et al (1959) use subadditive techniques to study the convergence of the
ratio~\(\frac{TSP_n}{\sqrt{n}} \longrightarrow \beta\) for some constant~\(\beta > 0,\) a.s.
as~\(n \rightarrow \infty.\) Another approach involves the study of concentration of~\(TSP_n\) around its mean
via concentration inequalities~(see Steele (1993)). Here we use the techniques
used in the proof of Theorem~\ref{tsp_thm} to obtain the following result.
\begin{Theorem}\label{var_tsp_thm} The variance
\begin{equation}\label{var_tsp_est_main}
\mathbb{E}\left(TSP_n - \mathbb{E} TSP_n\right)^2 \leq Cn^{2/3}
\end{equation}
for some constant~\(C > 0\) and for all~\(n \geq 1\) and so in particular,
\[\frac{1}{\sqrt{n}}\left(TSP_n - \mathbb{E} TSP_n\right)
\longrightarrow 0 \text{ in probability}\] as~\(n \rightarrow \infty.\)
Also there are positive constants~\(\{\theta_i\}_{1 \leq i \leq 3}\) such that
\begin{equation}\label{exp_tsp_u}
\theta_1 \sqrt{n} \leq \mathbb{E}TSP_n \leq 5 \sqrt{n},
\end{equation}
\begin{equation}\label{eq_tsp1_u}
\mathbb{P}\left(TSP_n \leq 5\sqrt{n} \right) =1
\end{equation}
and
\begin{equation}\label{eq_tsp2_u}
\mathbb{P}\left(TSP_n \geq \theta_2 \sqrt{n} \right) \geq 1 - \exp\left(-\frac{\theta_3 n}{\log{n}}\right)
\end{equation}
for all~\(n\) large.
\end{Theorem}

The paper is organized as follows. In Section~\ref{prelim}, we obtain preliminary estimates
needed for the proofs of main Theorems. In Section~\ref{pf_tsp}, we prove Theorem~\ref{tsp_thm}
and in Section~\ref{pf_tsp_var}, we prove Theorem~\ref{var_tsp_thm}.


\setcounter{equation}{0}
\renewcommand\theequation{\thesection.\arabic{equation}}
\section{Preliminary estimates}\label{prelim}
We first describe the strips method used throughout to find an upper bound for
the length of minimum length cycles.
\subsection*{Strips method}
Suppose there are~\(a \geq 3\) nodes~\(\{x_i\}_{1 \leq i\leq a}\) placed in a square~\(R\)
of side length~\(b.\) For~\(3 \leq j \leq a\) let~\(K(x_1,\ldots,x_j)\) be the
complete graph with vertex set~\(\{x_i\}_{1 \leq i \leq j}\)
and let~\({\cal C}_{j}\) be a spanning cycle of~\(K(x_1,\ldots,x_j)\) such that
\begin{equation}\label{min_cyc}
L({\cal C}_j) = \min_{\cal C} L({\cal C}) =: TSP(x_1,\ldots,x_j;R),
\end{equation}
where the minimum is taken over all spanning cycles of~\(K(x_1,\ldots,x_j)\) and~\(L({\cal C})\)
is the length of~\({\cal C}\) (see~(\ref{len_cyc_def})).

\begin{figure}
\centering
\begin{subfigure}{0.5\textwidth}
\centering
\includegraphics[width=3in, trim= 20 380 50 100, clip=true]{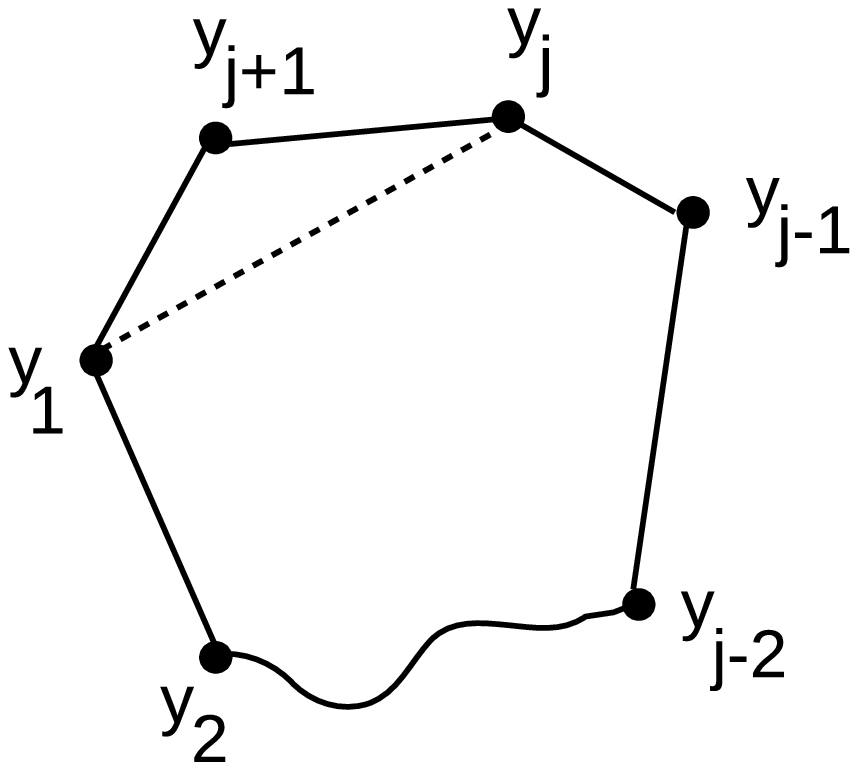}
   \caption{Removing the vertex~\(y_{j+1}\) from the cycle~\({\cal C}.\)}
\end{subfigure}
\begin{subfigure}{0.5\textwidth}
\centering
   \includegraphics[width=2in, trim= 20 180 50 100, clip=true]{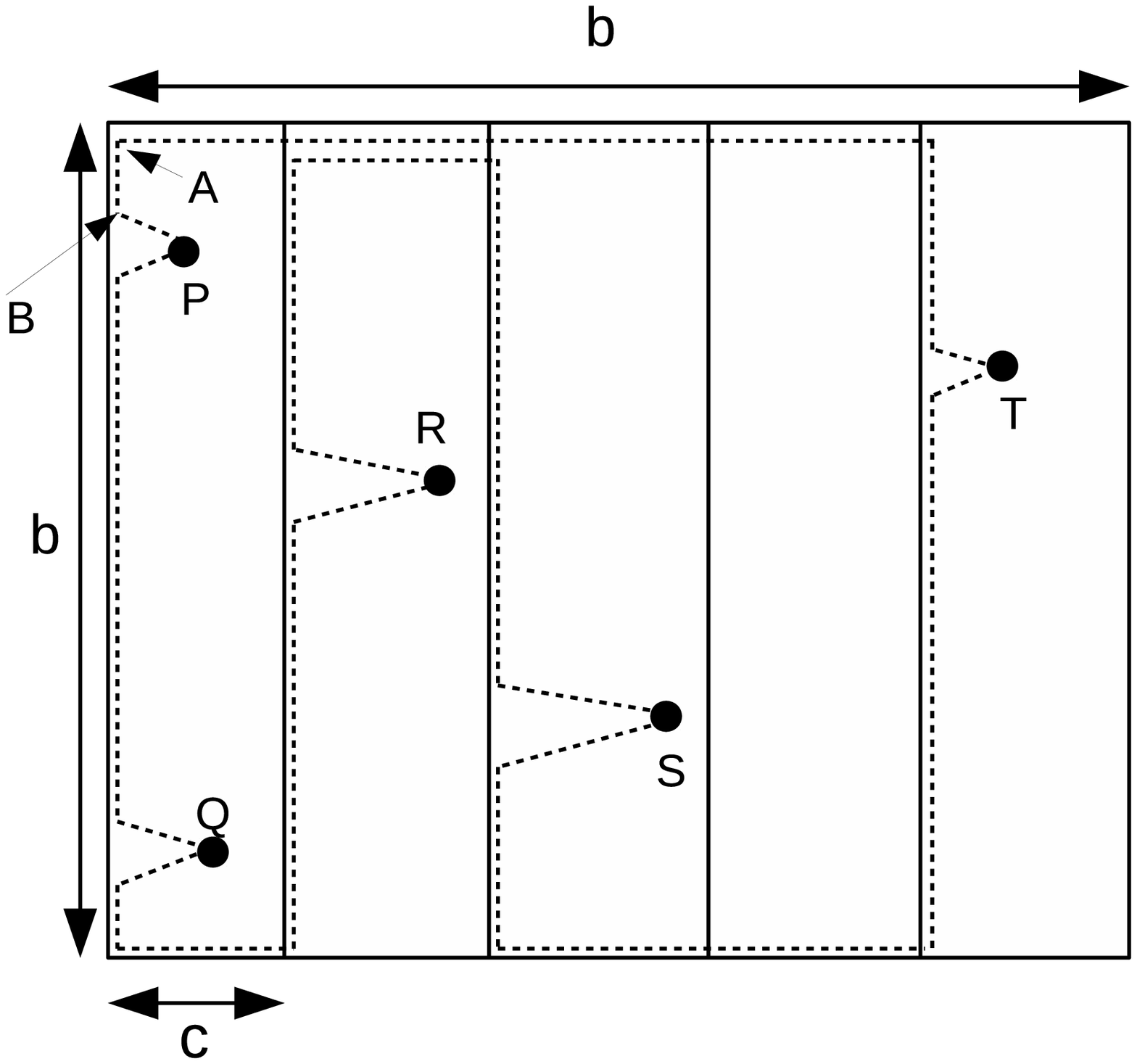}
   \caption{Estimating minimum length using strips counting.} 
  \end{subfigure}

\caption{\((a)\) Monotonicity in the TSP length. \((b)\) Estimating the length of the TSP using strips method.}
\label{stp_fig}
\end{figure}

For any~\(3 \leq j \leq a,\)
\begin{equation}\label{tsp_ab}
TSP(x_1,\ldots,x_j;R) \leq TSP(x_1,\ldots,x_a;R) \leq 5b \sqrt{a}.
\end{equation}
\emph{Proof of~(\ref{tsp_ab})}: The first estimate in~(\ref{tsp_ab}) is obtained
by monotonicity as follows. Let~\({\cal C} = (y_1,\ldots,y_{j+1},y_1)\) be any cycle
in~\(K(x_1,\ldots,x_a)\) with vertex set~\(\{y_i\}_{1 \leq i \leq j+1} = \{x_i\}_{1 \leq i \leq j+1}\)
and without loss of generality suppose that~\(y_{j+1} = x_{j+1}.\)
Recall that~\((y_j,y_{j+1})\) is the edge with~\(y_j\) and~\(y_{j+1}\)
as endvertices. Removing the edges~\((y_j,y_{j+1})\)
and~\((y_{j+1},y_1),\) and adding the edge~\((y_1,y_j)\)
we get a new cycle~\({\cal C}'\) with vertex set~\(\{x_i\}_{1 \leq i \leq j}\)
(see Figure~\ref{stp_fig}\((a)\)).

By triangle inequality, the lengths
\begin{equation}\label{trai1}
d(y_j,y_1) \leq d(y_j,y_{j+1}) + d(y_{j+1},y_1).
\end{equation}
and therefore the length~\(L({\cal C}')\) of~\({\cal C}'\)~(see (\ref{len_cyc_def}) for definition)
is
\[L({\cal C'}) = \sum_{i=1}^{j-1}d(y_i,y_{i+1}) + d(y_j,y_1)
\leq \sum_{i=1}^{j} d(y_i,y_{i+1}) + d(y_{j+1},y_1) = L({\cal C}).\]
But by definition~\(TSP(x_1,\ldots,x_j;b) \leq L({\cal C}')\)
and so~\(TSP(x_1,\ldots,x_j;b) \leq L({\cal C}).\)
Taking minimum over all cycles~\({\cal C}\)
with vertex set~\(\{x_i\}_{1 \leq i \leq j+1},\) we get
\[TSP(x_1,\ldots,x_j;R) \leq TSP(x_1,\ldots,x_{j+1};R).\]

For the second estimate in~(\ref{tsp_ab}), divide~\(R\) into vertical rectangles (strips) each of
size~\(c \times b\) so that the number of strips
is~\(\frac{b}{c}\) as shown in Figure~\ref{stp_fig}\((b).\) Here~\(a= 5\) and without loss of generality
suppose that~\(P = x_1, Q = x_2, R = x_3, S = x_4\) and~\(T = x_5.\)
The dotted line corresponds to a cycle
containing all the nodes~\(P,Q,R,S\) and~\(T.\) Starting from close to the top left corner
at point~\(A,\) we go vertically down and encounter the nodes~\(P,Q,R,S\) and~\(T\) in that order.
Each time we are close to a node, we ``reach" for the node by a slightly inclined line.
For example, the node~\(P\) is joined to the vertical dotted line~\(AB\) by the inclined
line~\(BP.\)

After the final node~\(T\) is encountered, we join it to the starting point~\(A\)
by inclined, vertical and horizontal lines as shown in Figure~\ref{stp_fig}. The cycle~\({\cal D}\) constructed
above consists of vertical, horizontal and inclined lines.
The number of strips is~\(\frac{b}{c}\) and the
sum of the lengths of the vertical lines in a particular strip is at most the height of the strip~\(b.\)
Therefore the total length of vertical lines in~\({\cal D}\) is at most~\(\frac{b}{c}b.\)

The total length of the horizontal lines in~\({\cal D}\) before encountering the final node~\(T\)
is at most~\(b.\) Since~\(T\) is joined to~\(A\) by a curve consisting of a
horizontal line, the total length of horizontal lines in~\({\cal D}\) is at most~\(2b.\)

Finally, each inclined line in~\({\cal D}\) has length at most~\(\frac{c}{\sqrt{2}},\) since the corresponding slope
is at most~\(45\) degrees. There are~\(a\) nodes and there are exactly two inclined lines
containing any particular node. Therefore the total length of the inclined lines in~\({\cal D}\)
is at most~\(a c \sqrt{2}.\)

Summarizing, the total length of edges in~\({\cal D}\) is at most~\(\frac{b^2}{c}  + ac\sqrt{2} + 2b.\)
By construction, the cycle~\({\cal D}\)
encounters the nodes~\(x_1,\ldots,x_a\) in that order and so applying triangle inequality
as before, the cycle~\({\cal C} = (x_1,x_2,\ldots,x_a,x_1)\)
with edges being the straight lines~\((x_1,x_{2}), (x_2,x_3),\ldots,(x_a,x_1),\)
has total length no more than the sum of length of edges in~\({\cal D}.\) Thus
\begin{equation}\label{temp_est_strip}
TSP(x_1,\ldots,x_a; R) \leq L({\cal C}) \leq \frac{b^2}{c}  + ac\sqrt{2} + 2b.
\end{equation}
Setting~\(c = \frac{b}{\sqrt{a}}\) in~(\ref{temp_est_strip}), we get
\[TSP(x_1,\ldots,x_a; R) \leq b\sqrt{a} + \sqrt{2}b\sqrt{a} + 2b \leq 5b\sqrt{a}, \]
since~\(a \geq 1.\)~\(\qed\)

\subsection*{Length of TSP within cities}
Recall from discussion prior to~(\ref{N_est}) that~\(n \geq 1\) nodes~\(\{X_k\}_{1 \leq k \leq n}\)
are distributed across the~\(r_n \times r_n\) squares~\(\{S_j\}_{1 \leq j \leq N}\)
according to a Binomial process with intensity~\(g_N\) as defined in~(\ref{gn_def}).
In this subsection, we obtain estimates for the length~\(T_l\)
of the minimum length cycle containing all the nodes of the square~\(S_l.\)

If~\(p_l\) denotes the probability that a node of~\(\{X_j\}\) occurs inside~\(S_l,\)
then
\begin{equation}\label{pl_def}
\frac{\eta_1}{N} \leq p_l := \frac{\int_{S_l} f(x) dx}{\int_{\cup_j S_j} f(x)  dx} \leq \frac{\eta_2}{N},
\end{equation}
where~\(\eta_1 = \frac{\epsilon_1}{\epsilon_2} \leq \frac{\epsilon_2}{\epsilon_1} = \eta_2\) (see~(\ref{f_eq})).
Therefore if
\begin{equation}\label{nl_def}
N_l = \sum_{i=1}^{n} \ind(X_i \in S_l)
\end{equation}
denotes the number of nodes of~\(\{X_j\}\) in the square~\(S_l,\) then~\(N_l\)
is Binomially distributed with parameters~\(n\) and~\(p_l;\) i.e., for any~\(1 \leq  k \leq n,\)
\begin{equation}\label{bin_dist}
\mathbb{P}(N_l = k) = B(k;n,p_{l}) := {n \choose k} p_l^{k} (1-p_l)^{n-k},
\end{equation}
where~\({n \choose k} = \frac{n!}{k!(n-k)!}\) is the Binomial coefficient.
Moreover,
\begin{equation}\label{nl_exp}
\frac{\eta_1 n}{N} \leq \mathbb{E} N_l = np_l \leq \frac{\eta_2 n}{N}
\end{equation}
by~(\ref{pl_def}).

Let~\(\{Y_j\}_{1 \leq j \leq N_l}\)
be the nodes of~\(\{X_j\}\) present in the square~\(S_l.\)
Formally, if~\(N_l = 0,\) set~\(\{Y_j\}_{1 \leq j \leq N_l} := \emptyset.\)
If~\(N_l \geq 1,\) define~\(N_l\) indices~\(j_1,\ldots,j_{N_l}\) as follows.
Let~\[j_1 = j_1(X_1,\ldots,X_n) := \min\{1 \leq k \leq n :  X_k \in S_l\}\]
be the least indexed node of~\(\{X_k\}\) present in~\(S_l.\)
Let~\[j_2 = \min\{j_1+1 \leq k \leq n : X_k \in S_l\}\] be the next least indexed node of~\(\{X_k\}\)
present in~\(S_l\) and so on. Set~\(Y_i = X_{j_i}\) for~\(1 \leq j \leq N_l.\)

Set~\(T_l = 0\) if~\(N_l \leq 2\) and if~\(N_l \geq 3\) set
\begin{equation}\label{tl_def}
T_l := TSP(Y_1,\ldots,Y_{N_l}; S_l)
\end{equation}
where~\(TSP(.;.)\) is as defined in~(\ref{min_cyc}). The following is the main lemma proved in this subsection.
\begin{Lemma} \label{tl_lemma}
If~\(M > 0\) is arbitrary and~(\ref{N_est}) holds, the following is true: There are positive constants~\(\{\delta_i\}_{1 \leq i \leq 3}\) such that
for all~\(n \geq 2\) and for any~\(1 \leq l \leq N,\)
\begin{equation}\label{del_tn_b}
\delta_1 r_n\sqrt{\frac{n}{N}}\leq \mathbb{E} T_l \leq \delta_2 r_n \sqrt{\frac{n}{N}} \;\;\; \text{ and }
\;\;\; \mathbb{E} T_l^2 \leq \delta_3 \left(r_n \sqrt{\frac{n}{N}}\right)^2.
\end{equation}
Moreover, if
\begin{equation}\label{ul_def}
U_l = U_l(n) := \left\{\frac{\eta_1 n}{2N} \leq N_l \leq \frac{2\eta_2 n}{N}\right\},
\end{equation}
where~\(\eta_1\) and~\(\eta_2\) are as in~(\ref{pl_def}),
then there are positive constants~\(\{\delta_i\}_{i=4,5}\) such that
for all~\(n \geq 2\) and for any~\(1 \leq l \leq N,\)
\begin{equation}\label{ul_est}
\mathbb{P}(U_l) \geq 1- \exp\left(-\delta_4 \frac{n}{N}\right)
\end{equation}
and
\begin{equation}\label{tl_bd_sm}
T_l\ind(U_l) \leq \delta_5 r_n\sqrt{\frac{n}{N}}.
\end{equation}
\end{Lemma}

To prove~ the above Lemma, we perform some preliminary computations.
We first derive bounds for the total number of squares~\(N.\) From~(\ref{N_est}) we have
that~\(r_n^2 \geq \frac{M\log{n}}{n}\) and since all the~\(r_n \times r_n\) squares~\(\{S_l\}_{1 \leq l \leq N}\)
are contained within the unit square~\(S,\) we also have~\(Nr_n^2 \leq 1\) and therefore~\(N \leq \frac{n}{M\log{n}}.\) Similarly from~(\ref{N_est})
we also have that~\(\frac{n}{N^2} \longrightarrow 0\) as~\(n \rightarrow \infty\) and so~\(N \geq \sqrt{n}\) for all~\(n\) large. Combining we get
\begin{equation}\label{n_N}
\sqrt{n} \leq N \leq \frac{n}{M\log{n}} \text{ and } \frac{n}{N} \geq M\log{n}
\end{equation}
for all~\(n\) large.

For~\(k \geq 2,\) let~\(D_l(k)\) be the expected minimum distance between the node~\(Y_k\)
and every other node in~\(S_l,\) given that there are~\(N_l = k\) nodes in~\(S_l;\) i.e.,
\begin{equation}\label{pf_exp1}
D(k) = D_{l}(k) := \mathbb{E}\left(d(Y_{k},\{Y_u\}_{1 \leq u \leq k-1}) | N_l = k\right),
\end{equation}
where~\(d(A,B) = \min_{x \in A, y \in B} d(x,y)\) is the minimum distance between finite sets~\(A\) and~\(B.\)
We have the following properties.\\
\((b1)\) For any~\(k \geq 2\) and~\(1 \leq l \leq N,\) the term
\begin{equation}\label{d_low}
D_l(k) \geq  \int_0^{\frac{r_n}{\sqrt{\delta}}} \left(1-\pi \eta_2 \left(\frac{r}{r_n}\right)^2\right)^{k-1} dr
\end{equation}
where~\(\eta_2 = \frac{\epsilon_2}{\epsilon_1}\) is as in~(\ref{pl_def}).\\\\
\((b2)\) There are positive constants~\(\gamma_i,1 \leq i \leq 3\) such that for any~\(k \geq 2\) and~\(1 \leq l \leq N,\) the minimum distance
\begin{equation}\label{x2_est}
\gamma_1 \frac{r_n}{\sqrt{k}} \leq D_l(k) \leq \gamma_2 \frac{r_n}{\sqrt{k}}
\text{ and } \mathbb{E}\left(d^2(Y_{k},\{Y_u\}_{1 \leq u \leq k-1}) | N_l = k\right) \leq \gamma_3 \frac{r_n^2}{k}.
\end{equation}
\emph{Proof of~\((b1)-(b2)\)}: Given~\(N_l = k,\) the nodes in~\(S_l\)
are independently distributed in~\(S_l\) with distribution~\(f;\) i.e.,
\begin{equation}\label{dl_zeq}
D_l(k) = \mathbb{E}d(Z_k,\{Z_j\}_{1 \leq j \leq k-1})
\end{equation}
where~\(\{Z_i\}_{1 \leq i \leq k}\) are i.i.d.\
with distribution
\begin{equation}\label{z_dist}
\mathbb{P}(Z_1 \in A) = \frac{\int_{A \cap S_l} f(x) dx}{\int_{S_l} f(x) dx}.
\end{equation}

Use Fubini's theorem and~(\ref{z_dist}) to write
\begin{equation}\label{dl_alt}
D_l(k) = \frac{1}{\int_{S_l} f(x) dx} \int_{S_l} \mathbb{E}d(x,\{Z_j\}_{1 \leq j \leq k-1})f(x)dx,
\end{equation}
where
\begin{equation}\label{edx1}
\mathbb{E}d(x,\{Z_j\}_{1 \leq j \leq k-1}) = \int_{0}^{\infty} \mathbb{P}(d(x,\{Z_j\}_{1 \leq j \leq k-1}) \geq r) dr.
\end{equation}
For any~\(x \in S_l,\) the minimum distance from~\(x\) to~\(\{Z_1,\ldots,Z_{k-1}\}\) is at least~\(r\) if and
only if~\(B(x,r) \cap S_l\) contains no point of~\(\{Z_j\}_{1 \leq j \leq k-1}.\) Here~\(B(x,r)\) is the
ball of radius~\(r\) centred at~\(x.\) Wherever the point~\(x \in S_l,\) the area of~\(B(x,r) \cap S_l\)
is at most~\(\pi r^2\) and so together with~(\ref{f_eq}), we then get
\begin{eqnarray}
\mathbb{P}(d(x,\{Z_j\}_{1 \leq j \leq k-1}) \geq r) &=& \left(1-\frac{\int_{B(x,r) \cap S_l} f(x)dx}{\int_{S_l} f(x)dx}\right)^{k-1}
\label{dx_min}\\
&\geq& \left(1-\pi \eta_2\frac{r^2}{r_n^2}\right)^{k-1}, \nonumber
\end{eqnarray}
where~\(\eta_2 =\frac{\epsilon_2}{\epsilon_1}\) is as in~(\ref{pl_def}).

To prove the lower bound for~\(D_l(k)\) in~(\ref{x2_est}) of~\((b2),\) fix \(k \geq 2\)  and use~(\ref{d_low}) to get that
\begin{equation}
D_l(k) \geq \int_{0}^{\frac{r_n}{\sqrt{\delta k}}}\left(1-\delta\left(\frac{r}{r_n}\right)^2\right)^{k-1} dr
\geq \int_{0}^{\frac{r_n}{\sqrt{\delta k}}}\left(1-\frac{1}{k}\right)^{k-1} dr \geq \frac{e^{-1} r_n}{\sqrt{\delta k}} \nonumber`
\end{equation}
for all~\(n\) large. The final estimate is obtained by using~\(\left(1-\frac{1}{r}\right)^{r-1} \geq e^{-1}\) for all~\(r \geq 2.\)

For the upper bound for~\(D_l(k)\) in~(\ref{x2_est}),  again use~(\ref{dx_min}) and the fact that~\(B(x,r) \cap S_l\) has area at least~\(\frac{\pi r^2}{4}\)
no matter where the position of~\(x,\) to get
\[\mathbb{P}(d(x,\{Z_j\}_{1 \leq j \leq k-1}) \geq r) \leq \left(1-\frac{\pi}{4\epsilon_1} \left(\frac{r}{r_n}\right)^2\right)^{k-1} \leq \exp\left(-\frac{\pi(k-1)}{4\epsilon_1 r_n^2} r^2\right)\] and so
\[D_l(k) \leq \int_{0}^{\infty} \exp\left(-\frac{\pi(k-1)}{4\epsilon_1 r_n^2} r^2\right) dr \leq \frac{C r_n}{\sqrt{k-1}} \leq \frac{2C r_n}{\sqrt{k}}\]
for all~\(k \geq 2\) and for some positive constant~\(C,\) not depending on~\(k\) or~\(l.\)

Finally for the second moment estimate in~(\ref{x2_est}), we argue analogous to~(\ref{pf_exp1}) and get that the
term~\(\mathbb{E}\left(d^2(Y_{k},\{Y_u\}_{1 \leq u \leq k-1}) | N_l = k\right) \) equals
\begin{equation}\label{pf_exp2}
\mathbb{E}d^2(Z_k,\{Z_j\}_{1 \leq j \leq k-1}) = \frac{1}{\int_{S_l} f(x) dx} \int_{S_l} \mathbb{E}d^2(x,\{Z_j\}_{1 \leq j \leq k-1}) f(x) dx
\end{equation}
where~\(\{Z_i\}_{1 \leq i \leq k}\) are i.i.d.\ with distribution as in~(\ref{z_dist}).
Arguing as in the previous paragraph we get
\begin{eqnarray}
\mathbb{E}(d^2(x,\{Z_j\}_{1 \leq j \leq k-1})) &=& \int r \mathbb{P}(d(x,\{Z_j\}_{1 \leq j \leq k-1}) \geq r ) dr \nonumber\\
&\leq& \int_{0}^{\infty} r\exp\left(-\frac{\pi(k-1)}{4\epsilon_1 r_n^2} r^2\right) dr \nonumber\\
&\leq& \frac{C r_n^2}{k} \label{tataq}
\end{eqnarray}
for some constant~\(C > 0\) not depending on~\(k\) or~\(x.\) Substituting~(\ref{tataq}) into~(\ref{pf_exp2}) gives
the desired bound for the second moment in~(\ref{x2_est}).~\(\qed\)


\emph{Proof of Lemma~\ref{tl_lemma}}: The proof of~(\ref{ul_est}) follows from standard Binomial estimates and the estimate for~\(\mathbb{E}N_l\) in~(\ref{nl_exp}). The proof of~(\ref{tl_bd_sm}) follows from the strips estimate~(\ref{tsp_ab}) with~\(a = \frac{2\eta_2 n}{N}\) and~\(b = r_n.\)


To prove the first estimate of~(\ref{del_tn_b}) assume~\(N_l \geq 3\) and
recall that~\(\{Y_u\}_{1 \leq u \leq N_l}\) are the nodes of the Binomial process in the square~\(S_l\) (see paragraph prior to~(\ref{pf_exp1})).
Let~\({\cal C}_l\) denote the minimum length cycle of length~\(T_l\)
containing the nodes~\(\{Y_u\}_{1 \leq u \leq N_l}.\)  If~\(l(Y_u,{\cal C}_l),1 \leq u \leq N_l\) is the sum of length of the two edges containing~\(Y_u\)
as an endvertex then \[l(Y_u,{\cal C}_l) \geq 2d(Y_u,\{Y_v\}_{v \neq u}),\] the minimum distance of~\(Y_u\) from all
the other nodes in~\(S_l\) as defined in~(\ref{pf_exp1}).

From~(\ref{len_cyc_def}),
\[T_l = L({\cal C}_l) = \frac{1}{2} \left(\sum_{u=1}^{N_l} l(Y_u,{\cal C}_l) \right) \geq \left(\sum_{u=1}^{N_l}d(Y_u,\{Y_v\}_{v  \neq u })\right)\]
and so
\begin{equation}
\mathbb{E} T_l  = \sum_{k \geq 3} \mathbb{E} T_l \ind(N_l = k) \geq \mathbb{E} \sum_{k \geq 3} \sum_{u=1}^{k} d(Y_u,\{Y_v\}_{v  \neq u })\ind(N_l = k). \label{ind_tl_b2}
\end{equation}
Recalling the definition of~\(D_l(k)\) in~(\ref{pf_exp1}) we further get
\begin{equation}
\mathbb{E} T_l  = \sum_{k \geq 3} \mathbb{P}(N_l = k) kD_l(k) \geq \sum_{\frac{\eta_1 n}{2N} \leq k \leq \frac{2\eta_2 n}{N}}\mathbb{P}(N_l = k) kD_l(k), \label{ind_tl_b}
\end{equation}
provided~\(n\) is large enough so that~\[\frac{\eta_1 n}{2N} \geq \frac{\eta_1}{2} M\log{n} \geq 3,\] the middle estimate being true because of~(\ref{n_N}). 

Using the estimate~\(D_l(k) \geq \frac{\gamma_1 r_n}{\sqrt{k}}\) (see~(\ref{x2_est})) in~(\ref{ind_tl_b}) we then get
\begin{eqnarray}
\mathbb{E} T_l  &\geq& \gamma_1 r_n\sum_{\frac{\eta_1 n}{2N} \leq k \leq \frac{2\eta_2 n}{N}}\mathbb{P}(N_l = k) \sqrt{k}  \nonumber\\
&\geq&  \gamma_1 r_n \sqrt{\frac{\eta_1 n}{2N}} \sum_{\frac{\eta_1 n}{2N} \leq k \leq \frac{2\eta_2 n}{N}}\mathbb{P}(N_l = k)  \nonumber\\
&\geq& \gamma_1 r_n \sqrt{\frac{\eta_1 n}{2N}}  \left(1-\exp\left(-C\frac{n}{N}\right)\right), \label{ind_t12_bn}
\end{eqnarray}
for some constant~\(C > 0,\) by~(\ref{ul_est}). Since~\(\frac{n}{N}  \longrightarrow \infty\) as~\(n \rightarrow \infty,\) (see~(\ref{n_N})), we get the lower bound for~\(\mathbb{E} T_l\) from~(\ref{ind_t12_bn}).


For the upper bound of~\(\mathbb{E}T_l\) in~(\ref{del_tn_b}), we argue as follows. Recall that~\(T_l = L({\cal C}_l)\) is the length
of the minimum length cycle~\({\cal C}_l\) containing all the~\(N_l\) nodes of~\(\{X_k\}\) in~\(S_l.\) If the number of nodes~\(N_l \leq \frac{2\eta_2 n}{N},\) then from~(\ref{tl_bd_sm}), we have that~\(T_l \leq C r_n \sqrt{\frac{n}{N}}\) for some constant~\(C > 0.\)
If~\(N_l \geq \frac{2\eta_2 n}{N},\) then use the fact that~\(T_l\) is bounded above by~\(N_l r_n \sqrt{2},\) since each edge in~\({\cal C}_l\)
has both endvertices in the~\(r_n \times r_n\) square~\(S_l\) and therefore has length at most~\(r_n\sqrt{2}.\) 
Thus
\begin{equation}\label{tij_est_temp_b}
\mathbb{E} T_l \leq C r_n \sqrt{\frac{n}{N}} + r_n\sqrt{2} \mathbb{E} \left(N_l\ind\left(N_l > \frac{2\eta_2 n}{N}\right)\right)
\leq C r_n \sqrt{\frac{n}{N}} + r_n\sqrt{2} \mathbb{E} (N_l\ind(U_l^c)),
\end{equation}
where~\(U_l\) is as defined in~(\ref{ul_def}).

Recall from discussion following~(\ref{nl_def}) that~\(N_l\) is Binomially distributed with parameters~\(n\) and~\(p_l\) and so by standard Binomial estimates~
\begin{equation}\label{enl2}
\mathbb{E}N_l^2 \leq C(np_l)^2 \leq \frac{Cn^2}{N^2}
\end{equation}
for some constant~\(C > 0,\) where the final estimate in~(\ref{enl2}) follows from the estimate for~\(p_l\) in~(\ref{pl_def}).
Using Cauchy-Schwarz inequality we therefore get
\begin{equation} \label{ul_2233}
\mathbb{E} N_l \ind(U_l^c) \leq \left(\mathbb{E}N_l^2\right)^{\frac{1}{2}} \left(\mathbb{P}(U_l^c)\right)^{\frac{1}{2}} \leq C_1\frac{n}{N}\exp\left(-C_2\frac{n}{N}\right) \leq \sqrt{\frac{n}{N}},
\end{equation}
for all~\(n\) large and for some positive constants~\(C_1,C_2.\) The middle inequality in~(\ref{ul_2233}) follows from~(\ref{ul_est})
and the final inequality in~(\ref{ul_2233}) is true since~\(\frac{n}{N} \longrightarrow \infty\) as~\(n \rightarrow \infty\) (see~(\ref{n_N})).
Substituting~(\ref{ul_2233}) into~(\ref{tij_est_temp_b}) gives the upper bound for~\(\mathbb{E}T_l\) in~(\ref{del_tn_b}). The proof of the bound for~\(\mathbb{E}T_l^2\) is analogous as above.~\(\qed\)

Define the covariance between~\(T_{l_1}\) and~\(T_{l_2}\) for distinct~\(l_1\) and~\(l_2\) as
\begin{equation}\label{cov_def}
cov(T_{l_1},T_{l_2}) = \mathbb{E}T_{l_1} T_{l_2} - \mathbb{E} T_{l_1} \mathbb{E}T_{l_2}.
\end{equation}
We need the following result for future use. Recall the definition of~\(\epsilon_1\) and~\(\epsilon_2\) in~(\ref{f_eq}).
\begin{Lemma}\label{cov_lemma}
There is a positive constant~\(M_0 = M_0(\epsilon_1,\epsilon_2)\) large so that the following holds if~(\ref{N_est})
is satisfied with~\(M > M_0:\) There are positive constants~\(C_1,C_2\) such that for all~\(n \geq 2\) and for any~\(1 \leq l_1 \neq l_2 \leq N,\)
\begin{equation}\label{cov_tl_est}
|cov(T_{l_1},T_{l_2})| \leq C_1 \left(\mathbb{E}T_{l_1} \mathbb{E}T_{l_2}\right)\frac{n}{N^2} \leq C_2 \frac{r_n^2 n^2}{N^3}.
\end{equation}
\end{Lemma}
To prove Lemma~\ref{cov_lemma}, we use Poissonization described in the next subsection.

\subsection*{Poissonization}
Recall from discussion prior to~(\ref{N_est}) that~\(n \geq 1\) nodes~\(\{X_k\}_{1 \leq k \leq n}\)
are distributed  across the~\(r_n \times r_n\) squares~\(\{S_j\}_{1 \leq j \leq N}\)
according to a Binomial process
with intensity~\(g_N(.)\) as defined in~(\ref{gn_def}).
Throughout, we use Poissonization as a tool
to obtain estimates for probabilities of events for the corresponding Binomial process.
We make precise the notions in this subsection.

Let~\({\cal P}\) be a Poisson process
on the squares~\(\cup_{j=1}^{N} S_j\) with intensity function~\(ng_N(.)\)
defined on the probability space~\((\Omega_0,{\cal F}_0, \mathbb{P}_0).\)
If~\(N^{(P)}_l\) be the number of nodes of~\({\cal P}\) present in the square~\(S_l, 1 \leq l \leq N,\)
then
\begin{equation}\label{poi_dist}
\mathbb{P}_0(N^{(P)}_l = k) = Poi(k;np_l) := e^{-np_l}\frac{(np_l)^{k}}{k!},
\end{equation}
where~\(p_l\) is as defined in~(\ref{pl_def}).
Moreover,
\begin{equation}\label{nl_exp_p}
\frac{\eta_1 n}{N} \leq \mathbb{E}_0 N^{(P)}_l = np_l \leq \frac{\eta_2 n}{N}
\end{equation}
by~(\ref{pl_def}).

Let~\(\{Y_j\}_{1 \leq j \leq N^{(P)}_l}\)
be the nodes of~\({\cal P}\) present in the square~\(S_l.\)
Analogous to~(\ref{tl_def}), set~\(T^{(P)}_l = 0\) if~\(N^{(P)}_l \leq 2\) and if~\(N^{(P)}_l \geq 3\) set
\begin{equation}\label{tl_def_p}
T^{(P)}_l := TSP(Y_1,\ldots,Y_{N^{(P)}_l}; S_l)
\end{equation}
where~\(TSP(.;.)\) is as defined in~(\ref{min_cyc}).
The following result is analogous to Lemma~\ref{tl_lemma}.
\begin{Lemma} \label{tl_lemma_poiss}
If~\(M > 0\) is arbitrary and~(\ref{N_est}) holds, then the following is true: There are positive constants~\(\{\delta_i\}_{1 \leq i \leq 3}\) such that
for all~\(n \geq 2\) and for any~\(1 \leq l \leq N,\)
\begin{equation}\label{del_tn}
\delta_1 r_n\sqrt{\frac{n}{N}}\leq \mathbb{E}_0 T^{(P)}_l \leq \delta_2 r_n \sqrt{\frac{n}{N}},\;\;\;\;\mathbb{E}_0 \left(T^{(P)}_l\right)^2 \leq \delta_3 \left(r_n \sqrt{\frac{n}{N}}\right)^2
\end{equation}
and
\begin{equation}\label{tn_prob_est}
\mathbb{P}_0\left(T^{(P)}_l \geq \delta_4 r_n \sqrt{\frac{n}{N}}\right) \geq \delta_5.
\end{equation}
\end{Lemma}
\emph{Proof of Lemma~\ref{tl_lemma_poiss}}: The proof of~(\ref{del_tn}) is analogous as in the Binomial case and proceeds as follows.
Define
\begin{equation}\label{ul_def_p}
U^{(P)}_l = U^{(P)}_l(n) := \left\{\frac{\eta_1 n}{2N} \leq N^{(P)}_l \leq \frac{2\eta_2 n}{N}\right\},
\end{equation}
where~\(\eta_1\) and~\(\eta_2\) are as in~(\ref{pl_def}). Analogous to~(\ref{ul_est}), the following bound
is obtained by standard Poisson distribution estimates: There is a positive constant~\(\gamma\) such that
for all~\(n \geq 2\) and for any~\(1 \leq l \leq N,\)
\begin{equation}\label{ul_est_p}
\mathbb{P}_0\left(U^{(P)}_l\right) \geq 1- \exp\left(-\gamma \frac{n}{N}\right).
\end{equation}

As in the Binomial case, given~\(N^{(P)}_l = k,\) the nodes of~\({\cal P}\) are i.i.d.\ distributed
according to distribution~(\ref{z_dist}). Therefore for~\(k \geq 2\)
we let~\[D^{(P)}_l(k) = \mathbb{E}_0\left(d(Y_k,\{Y_j\}_{1 \leq j \leq k-1})|N^{(P)}_l  = k\right)\]
and as in~(\ref{pf_exp1}) obtain that
\begin{equation}\label{dlp_k}
D^{(P)}_l(k) = \mathbb{E}d(Z_k,\{Z_j\}_{1 \leq j \leq k-1}) = D_l(k),
\end{equation}
where~\(D_l(k)\) is as defined in~(\ref{pf_exp1}), the random variables~\(\{Z_j\}_{1 \leq j \leq k}\)
are i.i.d.\ with distribution~(\ref{z_dist}) and the final equality in~(\ref{dlp_k}) is true
because of~(\ref{dl_zeq}). Consequently~\(D^{(P)}_l(k)\) also satisfies properties~\((b1)-(b2)\)
and the rest of the proof of~(\ref{del_tn}) is analogous to the Binomial case.

Finally, the estimate in~(\ref{tn_prob_est}) is obtained by using~(\ref{del_tn}) and the Paley-Zygmund inequality
\[\mathbb{P}_0\left(T^{(P)}_l \geq \lambda \mathbb{E}_0 T^{(P)}_l\right)
\geq (1-\lambda)^2\frac{(\mathbb{E}_0T^{(P)}_l)^2}{\mathbb{E}_0\left(T^{(P)}_l\right)^2}\] for~\(0 < \lambda <1.\)~\(\qed\)


We now use Poissonization
and obtain intermediate estimates needed to prove Lemma~\ref{cov_lemma}.
Recall from~(\ref{tl_def}) and~(\ref{tl_def_p}) that~\(T_l\) and~\(T^{(P)}_l\) are the lengths of the minimum length cycles containing
all the nodes in the~\(r_n \times r_n\) square~\(S_l, 1 \leq l \leq N\) in the Binomial and the Poisson process, respectively.
Recall the definition of~\(\epsilon_1\) and~\(\epsilon_2\) in~(\ref{f_eq}).
\begin{Lemma}\label{poi_bin_diff_lemm}
There is a positive constant~\(M_0 = M_0(\epsilon_1,\epsilon_2)\) large so that the following holds if~(\ref{N_est})
is satisfied with~\(M > M_0:\) There are positive constants~\(C_0,C_1\) and~\(C_2\) such that for all~\(n \geq C_0\) and for any~\(1 \leq l \leq N,\)\\
\begin{equation}\label{tl_diff}
|\mathbb{E}T_l - \mathbb{E}_0 T^{(P)}_l| \leq C_1\left(\mathbb{E}T_l\right)\left(\frac{n}{N^2}\right)
\leq C_2\left(\frac{r_n n^{3/2}}{N^{5/2}}\right).
\end{equation}
Moreover, for any~\(1 \leq l_1 \neq l_2 \leq N\)
\begin{equation}\label{tl_diff2}
|\mathbb{E}(T_{l_1}T_{l_2}) - \mathbb{E}_0 (T^{(P)}_{l_1} T^{(P)}_{l_2})| \leq C_1\left(\mathbb{E}T_{l_1}\mathbb{E}T_{l_2}\right)\left(\frac{n}{N^2}\right)
\leq C_2\left(\frac{r_n^2 n^2}{N^3}\right).
\end{equation}
\end{Lemma}

To prove Lemma~\ref{poi_bin_diff_lemm}, we need estimates on the difference between Binomial and Poisson distributions.
For~\(k,l \geq 1\) recall the Binomial distribution~\(B(k; n,p_l)\)
and the Poisson distribution~\(Poi(k; np_l)\)
as defined in~(\ref{bin_dist}) and~(\ref{poi_dist}), respectively.
For~\(k_1,k_2,l_1,l_2 \geq 1,\) let
\begin{equation}\label{bin_dist2}
B(k_1,k_2;n,p_{l_1},p_{l_2}) := {n \choose k_1,k_2} p_{l_1}^{k_1} p_{l_2}^{k_2}(1-p_{l_1} - p_{l_2})^{n-k_1-k_2},
\end{equation}
where~\({n \choose k_1,k_2} = \frac{n!}{k_1!k_2!(n-k_1-k_2)!}.\)
We have the following properties.\\
\((c1)\) There is a constant~\(C > 0\) such that for all~\(n \geq 3,\)~\(1 \leq l \leq N\) and~\(\frac{\eta_1 n}{2N} \leq k \leq \frac{2\eta_2 n}{N},\)
\begin{equation}\label{b_poi_diff}
|B(k;n,p_{l}) - Poi(k;np_{l})| \leq Poi(k;np_{l}) \left(1+ \frac{C n}{N^2}\right).
\end{equation}
\((c2)\) There is a constant~\(C > 0\) such that for all~\(n \geq 3,\)
and for any~\(1 \leq l_1,l_2 \leq N\) and~\(\frac{\eta_1 n}{2N} \leq k_1,k_2 \leq \frac{2\eta_2 n}{N},\)
\begin{eqnarray}
&&|B(k_1,k_2;n,p_{l_1},p_{l_2}) - Poi(k_1;np_{l_1})Poi(k_2; np_{l_2})| \nonumber\\
&&\;\;\;\;\leq Poi(k_1;np_{l_1}) Poi(k_2;np_{l_2})\left(1+ \frac{C n}{N^2}\right). \label{b_poi_diff2}
\end{eqnarray}
\emph{Proof of~\((c1)-(c2)\)}:
To prove~(\ref{b_poi_diff}) in~\((c1),\) we write~\(p_{l} =p\) for simplicity.
Use~\({n \choose k} \leq \frac{n^{k}}{k!}\) and~\(1-x \leq e^{-x}\) for~\(0 < x < 1\) to get
\[{n \choose k} p^{k} (1-p)^{n-k} \leq \frac{(np)^{k}}{k!}e^{-p(n-k)} = Poi(k;np)e^{kp}.\]
Using~(\ref{pl_def}) and the fact that~\(k \leq \frac{2\eta_2}{N}\) we get \[e^{kp} \leq \exp\left(\frac{k\eta_2 n}{N}\right)
\leq \exp\left(2\eta_2 \frac{n}{N^2}\right)\]
and since~
\begin{equation}\label{ex_est_small_x}
e^{x} = 1+x+\sum \frac{x^{k}}{k!} \leq 1+x+ \sum_{k \geq 2}x^{k} \leq 1+2x
\end{equation}
for all~\(x\) small, we get~\(e^{kp} \leq 1 + \frac{4\eta_2 n}{N^2},\) proving the upper bound in~(\ref{b_poi_diff}).

To obtain a lower bound, we use the estimate
\begin{equation}\label{exp_low}
1-x \geq e^{-x-x^2}
\end{equation}
for all~\(0 < x< \frac{1}{2}.\) To prove~(\ref{exp_low}),
write \(\log(1-x) = -x - R(x)\) where~\[R(x) = \sum_{k \geq 2}\frac{x^{k}}{k} \leq \frac{1}{2}\sum_{k \geq 2}x^{k} =  \frac{x^2}{2(1-x)} \leq x^2\]
since~\(x < \frac{1}{2}.\)
Use~\({n \choose k} \geq \frac{(n-k)^{k}}{k!}\) and~(\ref{exp_low}) to get
\begin{equation}\label{bin2}
B(k;n,p) \geq \frac{1}{k!}(n-k)^{k} p^{k}e^{-p(n-k) - p^2(n-k)} = Poi(k;np)\left(1-\frac{k}{n}\right)^{k}e^{kp-(n-k)p^2}
\end{equation}
As before, using the fact that~\(\frac{\eta_1 n}{2N} \leq k \leq \frac{2\eta_2 n}{N}\) we get
\begin{equation}\label{temp_est2}
\left(1-\frac{k}{n}\right)^{k} \geq 1-\frac{k^2}{n} \geq 1-\frac{4\eta_2^2 n}{N^2}
\end{equation}
and~using~(\ref{pl_def}) we get
\begin{equation}
kp - (n-k)p^2 \geq kp - np^2 \geq \frac{\eta_1 n}{2N} \frac{\eta_1}{N}  - n\left(\frac{\eta_2 }{N}\right)^2 = -\eta \frac{n}{N^2}\label{temp_estt}
\end{equation}
where~\(\eta  = \eta^2_2 - \frac{\eta_1^2}{4}> 0,\) since~\(\epsilon_1 \leq  \epsilon_2\) and so~\(\eta_1  = \frac{\epsilon_1}{\epsilon_2} \leq \frac{\epsilon_2}{\epsilon_1} = \eta_2.\)
Using~(\ref{temp_est2}) and~(\ref{temp_estt}) into~(\ref{bin2}) gives
\begin{eqnarray}
B(k;n,p) &\geq& Poi(k;np)\left(1-\frac{\eta_1^2}{4} \frac{n}{N^2}\right) \exp\left(-\eta \frac{n}{N^2}\right) \nonumber\\
&\geq& Poi(k;np)\left(1-\frac{\eta_1^2}{4} \frac{n}{N^2}\right) \left(1-\eta \frac{n}{N^2}\right), \nonumber
\end{eqnarray}
since~\(e^{-x} \geq 1-x\) for~\(0 < x < 1.\) This proves~(\ref{b_poi_diff}).

To prove~(\ref{b_poi_diff2}), write~\(p_{l_1} =p_1, p_{l_2} = p_2\) and~\(B_{12} = B(k_1,k_2; n,p_1,p_2)\) for simplicity.
Use
\begin{equation}\label{eq_n_k1k2}
{n \choose k_1,k_2}  = \frac{1}{k_1!k_2!}n(n-1)\ldots (n-k_1-k_2+1) \leq \frac{n^{k_1+k_2}}{k_1!k_2!}
\end{equation}
to get
\begin{equation}\label{b12_est1}
B_{12} \leq \frac{(np_1)^{k_1}}{k_1!} \frac{(np_2)^{k_2}}{k_2!} e^{-(p_1+p_2)n}e^{(p_1+p_2)(k_1+k_2)}.
\end{equation}
Using~(\ref{pl_def}), we get~\(p_1+p_2 \leq \frac{2\eta_2 }{N}\) and since~\(k_1,k_2 \leq \frac{2\eta_2 n}{N}\)
we get using~(\ref{ex_est_small_x}) that
\begin{equation}\label{ep1p2}
e^{(p_1+p_2)(k_1+k_2)} \leq \exp\left(\frac{4\eta_2^2 n}{N^2}\right) \leq 1+\frac{8\eta_2^2 n}{N^2}
\end{equation}
for all~\(n\) large, since~\(\frac{n}{N^2} \longrightarrow 0\) as~\(n \rightarrow \infty\) (see~(\ref{N_est})).
Substituting~(\ref{ep1p2}) into~(\ref{b12_est1}), we get the upper bound for~\(B_{12}\) in~(\ref{b_poi_diff2}).

For the lower bound for~\(B_{12}\) again use~(\ref{eq_n_k1k2}) to get
\begin{equation}\nonumber
{n \choose k_1,k_2}  \geq \frac{1}{k_1!k_2!} (n-k_1-k_2)^{k_1+k_2} = \frac{n^{k_1+k_2}}{k_1!k_2!} \left(1-\frac{k_1+k_2}{n}\right)^{k_1+k_2}.
\end{equation}
Using~\((1-x)^{r} \geq 1-rx\) for~\(r,x > 0\) we further get
\begin{equation}\label{eq_n_k1k22}
{n \choose k_1,k_2} \geq \frac{n^{k_1+k_2}}{k_1!k_2!}  \left(1-\frac{(k_1+k_2)^2}{n}\right) \geq \frac{n^{k_1+k_2}}{k_1!k_2!}  \left(1-\frac{4\eta_2^2n}{N^2}\right)
\end{equation}
since~\(k_1,k_2 \leq \frac{2\eta_2n}{N}.\)
Substituting~(\ref{eq_n_k1k22}) into~(\ref{bin_dist2}) we get
\begin{equation}\label{b12_k12}
B_{12} \geq \frac{(np_1)^{k_1}}{k_1!} \frac{(np_2)^{k_2}}{k_2!} (1-p_1-p_2)^{n-k_1-k_2}\left(1-\frac{4\eta_2^2n}{N^2}\right).
\end{equation}

To evaluate~\((1-p_1-p_2)^{n-k_1-k_2},\) we use the estimate~(\ref{exp_low}) which is applicable since from~(\ref{pl_def}), we have~\[p_1+p_2 \leq \frac{2\eta_2}{N} \leq \frac{2\eta_2}{\sqrt{n}} \longrightarrow 0\] as~\(n \rightarrow \infty\) (see~(\ref{n_N})).
Using~(\ref{exp_low}), we get
\begin{equation}\label{p12_k12_est1}
(1-p_1-p_2)^{n-k_1-k_2} \geq e^{-(p_1+p_2)(n-k_1-k_2) - (p_1+p_2)^2(n-k_1-k_2)} = e^{-np_1} e^{-np_2} e^{I_{1}-I_2},
\end{equation}
where
\begin{equation}\label{i1_k12}
I_{1}= (p_1+p_2)(k_1+k_2) \geq 0
\end{equation}
and
\begin{equation}\label{i2_k12}
I_2 = (p_1+p_2)^2 (n-k_1-k_2) \leq  n(p_1+p_2)^2 \leq \frac{\eta_2^2n}{N^2}
\end{equation}
for some constant~\(C_1 > 0.\)   The final estimate in~(\ref{i2_k12}) follows from the fact that~\(p_1+p_2 \leq \frac{2\eta_2 n}{N}\) (see~(\ref{pl_def})). Using~\(e^{-x} \geq 1-x\) we get
\begin{equation}\label{i1_i2_k12}
e^{I_1-I_2} \geq e^{-I_2} \geq 1 - \frac{\eta_2^2 n}{N^2}
\end{equation}
and substituting~(\ref{i1_i2_k12}) into~(\ref{p12_k12_est1}), we
\begin{equation}\label{p12_k12_est2}
(1-p_1-p_2)^{n-k_1-k_2} \geq  e^{-np_1} e^{-np_2} \left(1-\frac{\eta_2^2 n}{N^2}\right).
\end{equation}
Using~(\ref{p12_k12_est2}) in~(\ref{b12_k12}), we get the lower bound for~\(B_{12}\) in~(\ref{b_poi_diff2}).~\(\qed\)


Using properties~\((c1)-(c2)\) we prove Lemma~\ref{poi_bin_diff_lemm}.\\
\emph{Proof of~(\ref{tl_diff}) in Lemma~\ref{poi_bin_diff_lemm}}: Recall from~(\ref{nl_def}) that~\(N_l\) is the number of nodes
of the Binomial process~\(\{X_k\}\) in the square~\(S_l\) and let~\(U_l\) be the event as defined in~(\ref{ul_def}). Write
\begin{equation}\label{etl_1}
\mathbb{E}T_l = I_1 + I_2
\end{equation}
where~\[I_1 = \mathbb{E} T_l \ind(U_l) = \sum_{\frac{\eta_1 n}{2N} \leq k \leq \frac{2\eta_2 n}{N}} \mathbb{E} T_l \ind(N_l = k),\] \(I_2 = \mathbb{E}T_l \ind(U^c_l)\) and~\(\eta_1,\eta_2\) are as in~(\ref{pl_def}).

Similarly
\begin{equation}\label{et_poi1}
\mathbb{E}_0T^{(P)}_l = I^{(P)}_1 + I^{(P)}_2
\end{equation}
where
\[I^{(P)}_1 = \mathbb{E}_0(T^{(P)}_l \ind(U^{(P)}_l)),\] \(I^{(P)}_2 = \mathbb{E}_0(T^{(P)}_l \ind(U^{(P)}_l)^c),\)
\(U^{(P)}_l = \left\{\frac{\eta_1 n}{2N} \leq N^{(P)}_l \leq \frac{2\eta_2 n}{N}\right\}\) is as defined in~(\ref{ul_def_p})
and~\(N^{(P)}_l\) is the number of nodes of the Poisson process~\({\cal P}\) inside the square~\(S_l\) (see discussion prior to~(\ref{poi_dist})).

From~(\ref{etl_1}) and~(\ref{et_poi1}), we therefore get
\begin{equation}\label{tl_diff33}
|\mathbb{E}T_l - \mathbb{E}_0T^{(P)}_l| \leq |I_1 - I^{(P)}_1| + I_2 + I^{(P)}_2.
\end{equation}
The remainder terms~\(I_2\) and~\(I^{(P)}_2\) satisfy
\begin{equation}\label{i2_est_fin2}
\max(I_2,I^{(P)}_2) \leq C (\mathbb{E} T_l) \frac{n}{N^2}
\end{equation}
for some constant~\(C > 0.\) We prove~(\ref{i2_est_fin2}) for~\(I_2\)
and an analogous proof holds for~\(I^{(P)}_2.\)
Indeed, every edge in the minimum length cycle~\({\cal C}_l\)
containing all the nodes in the~\(r_n \times r_n\) square~\(S_l\)
has both endvertices within~\(S_l\) and so has length at most~\(r_n\sqrt{2}.\)
Since there are~\(N_l\) nodes in the square~\(S_l,\) we must have~\(T_l \leq N_l r_n\sqrt{2}\)
and so
\begin{equation}\label{i2_etr22}
I_2  = \mathbb{E}T_l \ind(U_l^c) \leq r_n \sqrt{2} \mathbb{E}N_l\ind(U_l^c).
\end{equation}
Using the third expression in~(\ref{ul_2233}) to estimate~\(\mathbb{E}N_l \ind(U_l^c)\) we get
\begin{equation}\label{i2_etr2}
I_2  \leq C_1 r_n\sqrt{2} \frac{n}{N}\exp\left(-C_2 \frac{n}{N}\right) = C_1\sqrt{2} \left(r_n \sqrt{\frac{n}{N}}\right) \left(\sqrt{\frac{n}{N}}\exp\left(-C_2 \frac{n}{N}\right)\right)
\end{equation}
for some constants~\(C_1,C_2 > 0.\) From the lower bound in~(\ref{del_tn_b}) we have~\(\mathbb{E} T_l \geq C_3 r_n \sqrt{\frac{n}{N}}\)
and so
\begin{eqnarray}
I_2 &\leq& C_4 \left(\mathbb{E}T_l\right)\left(\sqrt{\frac{n}{N}}\exp\left(-C_2 \frac{n}{N}\right)\right) \nonumber\\
&=& C_4 \left(\mathbb{E}T_l\right)
\left(\frac{n}{N^2}\right) \left( \left(\frac{N^3}{n}\right)\exp\left(-\frac{C_2 n}{2N}\right)\right)^\frac{1}{2} \label{i2_etr33}
\end{eqnarray}

Using the upper bound~\(N \leq \frac{n}{M\log{n}}\) from~(\ref{n_N}), we have
\begin{equation}\label{m_bds2}
\left(\frac{N^3}{n}\right)\exp\left(-\frac{C_2 n}{2N}\right) \leq \frac{n^2}{M^3(\log{n})^3} \exp\left(-\frac{C_2M}{2}\log{n}\right) \leq 1
\end{equation}
for all~\(n\) large, provided~\(M >0\) large. Fixing such an~\(M\) and using~(\ref{m_bds2}) in~(\ref{i2_etr33}), we get~(\ref{i2_est_fin2}).


To estimate the difference~\(I_1 - I^{(P)}_1\) in~(\ref{tl_diff33}),
recall that given~\(N_l = k,\) the nodes in~\(S_l\) are independently distributed in~\(S_l\) with distribution~\(\frac{f(.)}{\int_{S_l} f(x)dx}\) (see~(\ref{z_dist})) and so
\begin{equation}
I_1 = \sum_{\frac{\eta_1 n}{2N} \leq k \leq \frac{2\eta_2 n}{N}}\mathbb{P}(N_l = k)\mathbb{E}(T_l| N_l =k)  = \sum_{\frac{\eta_1 n}{2N} \leq k \leq \frac{2\eta_2 n}{N}} B(k; n,p_l) \Delta(k,q_l) \label{i1_est1}
\end{equation}
where~\(B(k;n,p_{l})\) is the Binomial probability distribution as defined in~(\ref{bin_dist}),~\(q_l = \int_{S_l} f(x) dx,\)
\begin{equation}\label{del_p_def}
\Delta(k,q_l) = \mathbb{E}(T_l|N_l = k) = \int_{S_l} TSP(z_1,\ldots,z_k;S_l) \frac{f(z_1)}{q_l}\ldots \frac{f(z_k)}{q_l} dz_1\ldots dz_k
\end{equation}
and~\(TSP(z_1,\ldots,z_k;S_l)\) is the minimum length of a cycle containing all the nodes~\(z_1,\ldots,z_k \in S_l\) (see~(\ref{min_cyc})).

Similarly, as argued in~(\ref{dlp_k}), given~\(N_l^{(P)} = k,\) the nodes of the Poisson process~\({\cal P}\) are
also distributed in~\(S_l\) according to distribution~\(\frac{f(.)}{\int_{S_l} f(x)dx}.\) Therefore
\[\mathbb{E}(T^{(P)}_l | N^{(P)}_l = k) = \Delta(k,q_l)\] as defined in~(\ref{del_p_def}) and so
\begin{equation}\label{ip_est11}
I^{(P)}_1 = \sum_{\frac{\eta_1 n}{2N} \leq k \leq \frac{2\eta_2 n}{N}}\Delta(k,q_l) Poi(k;np_l),
\end{equation}
where~\(Poi(k;np_l)\) is the Poisson distribution as defined in~(\ref{poi_dist}).
From~(\ref{i1_est1}) and~(\ref{ip_est11}), we therefore get
\begin{equation}\label{i1p_diff}
|I_1 - I^{(P)}_1| \leq \sum_{\frac{\eta_1 n}{2N} \leq k \leq \frac{2\eta_2 n}{N}} \Delta(k,q_l) |B(k;n,p_l) - Poi(k;np_l)|.
\end{equation}

Using estimate~(\ref{b_poi_diff}) of property~\((c1)\) to approximate the Binomial distribution with the Poisson distribution, we get
\begin{eqnarray}
|I_1-I^{(P)}_1| &\leq& C_1 \left(\sum_{\frac{\eta_1 n}{2N} \leq k \leq \frac{2\eta_2 n}{N}} Poi(k;np_l) \Delta(k,q_l) \right)\frac{n}{N^2} \nonumber\\
&\leq& C_1 \left(\sum_{k \geq 0}Poi(k;np_l) \Delta(k,q_l)\right) \frac{n}{N^2} \nonumber\\
&=& C_1 \left(\mathbb{E}_0(T^{(P)}_l)\right) \frac{n}{N^2} \label{i1_temp_estt}
\end{eqnarray}
for some constant~\(C_1 >0.\) Finally, from~(\ref{del_tn_b}) and~(\ref{del_tn}), we obtain that both~\(\mathbb{E}_0(T^{(P)}_l)\) and~\(\mathbb{E}T_l\) are bounded above and below by constant multiples of~\(r_n \sqrt{\frac{n}{N}}\) and so~\(\mathbb{E}_0(T^{(P)}_l) \leq C_2 \mathbb{E} T_l\) for some constant~\(C_2 > 0\) and from~(\ref{i1_temp_estt}), we therefore get
\begin{equation}\label{i1_fin}
|I_1-I^{(P)}_1| \leq C_3 \left(\mathbb{E} T_l\right) \frac{n}{N^2}
\end{equation}
for some constant~\(C_3 >0.\)
Substituting~(\ref{i1_fin}) and~(\ref{i2_est_fin2}) into~(\ref{tl_diff33}) gives
\[|\mathbb{E}T_l- \mathbb{E}_0T^{(P)}_l| \leq C_4 \left(\mathbb{E} T_l\right) \frac{n}{N^2} \leq C_5 \left(\frac{r_nn^{3/2}}{N^{5/2}}\right),\] for some positive constants~\(C_4,C_5,\) again using the upper bound for~\(\mathbb{E}T_l\) from~(\ref{del_tn_b}). This proves~(\ref{tl_diff}).~\(\qed\)

\emph{Proof of~(\ref{tl_diff2}) of Lemma~\ref{poi_bin_diff_lemm}}: Recall the definition of~\(U_l\) in~(\ref{ul_def}) and write
\begin{equation}\label{t12_1}
\mathbb{E} T_{l_1} T_{l_2} = J_1 + J_2,
\end{equation}
where~\(J_1 = \mathbb{E} T_{l_1} T_{l_2} \ind( U_{l_1} \cap U_{l_2})\) and~\(J_2 = \mathbb{E}T_{l_1} T_{l_2} \ind(U^c_{l_1} \cup U^{c}_{l_2}).\)
Similarly, for the Poisson case let~\(U^{(P)}_l\) be the event defined in~(\ref{ul_def_p}) and write
\begin{equation}\label{t12_p1}
\mathbb{E}_0 T^{(P)}_{l_1} T^{(P)}_{l_2} = J^{(P)}_1 + J^{(P)}_2,
\end{equation}
where~\(J^{(P)}_1 = \mathbb{E}_0 T^{(P)}_{l_1} T^{(P)}_{l_2} \ind( U^{(P)}_{l_1} \cap U^{(P)}_{l_2})\) and~\(J^{(P)}_2 = \mathbb{E}_0T^{(P)}_{l_1} T^{(P)}_{l_2} \ind(U^{(P)}_{l_1} \cup U^{(P)}_{l_2})^c.\)

From~(\ref{t12_1}) and~(\ref{t12_p1}), we get
\begin{equation}\label{r12_diff}
|\mathbb{E}T_{l_1}R_{l_2} - \mathbb{E}_0 T^{(P)}_{l_1} T^{(P)}_{l_2}| \leq |J_1 - J_1^{(P)}| + J_2 + J_2^{(P)}.
\end{equation}
The remainder terms~\(J_2\) and~\(J^{(P)}_2\) satisfy
\begin{equation}\label{i12_est_fin2}
\max(J_2,J^{(P)}_2) \leq C_1 (\mathbb{E} T_{l_1}\mathbb{E}T_{l_2}) \frac{n}{N^2} \leq C_2 \left(\frac{r_n^2 n^2}{N^3}\right)
\end{equation}
for some constants~\(C_1,C_2 > 0.\) We prove~(\ref{i12_est_fin2}) for~\(J_2\)
and an analogous proof holds for~\(J^{(P)}_2.\)
As argued in the proof of~(\ref{i2_est_fin2}), every
one of the~\(N_{l_1}\) edges in the minimum length cycle~\({\cal C}_{l_1}\) of length~\(T_{l_1}\)
has both endvertices within~\(S_{l_1}\) and so has length at most~\(r_n\sqrt{2}.\)
Therefore
\begin{equation}\label{i12_etr22}
J_2  = \mathbb{E}T_{l_1}T_{l_2} \ind(U_{l_1}^c \cup U_{l_2}^c) \leq \left(r_n \sqrt{2}\right)^2 \mathbb{E}N_{l_1}N_{l_2}\ind(U^c_{l_1} \cup U^{c}_{l_2}).
\end{equation}
Using Cauchy-Schwarz inequality,
\begin{equation}\label{nl12_1}
\mathbb{E}N_{l_1}N_{l_2}\ind(U^c_{l_1} \cup U^{c}_{l_2}) \leq \left(\mathbb{E}N^2_{l_1}N^2_{l_2}\right)^{\frac{1}{2}}
\mathbb{P}\left(U^c_{l_1} \cup U_{l_2}^c\right)^{\frac{1}{2}}
\end{equation}
and using the estimate~(\ref{ul_est}), we have
\begin{equation}\label{ul12}
\mathbb{P}\left(U^c_{l_1} \cup U_{l_2}^c\right) \leq \mathbb{P}\left(U^c_{l_1}\right) +\mathbb{P}\left(U_{l_2}^c\right) \leq 2\exp\left(-4C \frac{n}{N}\right)
\end{equation} for some constant~\(C > 0\) and for all~\(n\) large.

To evaluate~\(\mathbb{E} N^2_{l_1} N^2_{l_2},\) use~\(ab \leq \frac{a^2+b^2}{2}\) to get
\begin{equation}\label{nl12_2}
\mathbb{E} N^2_{l_1} N^2_{l_2} \leq \frac{1}{2}\left(\mathbb{E}N_{l_1}^4 + \mathbb{E} N_{l_2}^4\right)
\end{equation}
and use the fact that the term~\(N_l\) is Binomially distributed with parameters~\(n\) and~\(p_l,\) where~\(p_l \leq  \frac{\eta_2}{N}\) (see~(\ref{pl_def}))
and~\(\eta_2\) does not depend on~\(l\) or~\(n.\) Therefore
~\[\mathbb{E}N_l^4 \leq C_1 (np_l)^4 \leq C_2 \left(\frac{n}{N}\right)^4\] for some constants~\(C_1,C_2\) not depending on~\(l\) or~\(n\)
and so from~(\ref{nl12_2}) we get
\begin{equation}\label{nl12_3}
\mathbb{E} N^2_{l_1} N^2_{l_2} \leq C_3 \left(\frac{n}{N}\right)^4.
\end{equation}
Using~(\ref{nl12_3}) and~(\ref{ul12}) in~(\ref{nl12_1})
we get
\begin{equation}\label{nl12_4}
\mathbb{E}N_{l_1}N_{l_2}\ind(U^c_{l_1} \cup U^{c}_{l_2}) \leq C_4\left(\frac{n}{N}\right)^2 \exp\left(-2C\frac{n}{N}\right).
\end{equation}

Substituting~(\ref{nl12_4}) into~(\ref{i12_etr22}) gives~(\ref{tl_diff2}).~\(\qed\)
\begin{equation}\label{i12_etr33}
J_2  \leq C_5 r_n^2 \left(\frac{n}{N}\right)^2 \exp\left(-2C\frac{n}{N}\right) = C_5 \left(\frac{r_n^2n^2}{N^3}\right) N\exp\left(-2C\frac{n}{N}\right).
\end{equation}
Since~\(N \leq \frac{n}{M\log{n}}\) (see~(\ref{n_N})) we have that
\[N\exp\left(-2C\frac{n}{N}\right) \leq \frac{n}{M\log{n}} \exp\left(-2CM\log{n}\right) \leq 1\] for all~\(n\) large provided~\(M > 0\) is large.
Fixing such an~\(M,\) we get~(\ref{i12_est_fin2}).

To evaluate the difference~\(J_1-J^{(P)}_1,\) recall from discussion prior to~(\ref{i1_est1}) that given~\(N_l = k,\)
the nodes of the Binomial process are distributed in the square~\(S_l\) with distribution~(\ref{z_dist}). Similarly,
given~\(N_l^{(P)} = k,\) the nodes of the Poisson process are also distributed according to~(\ref{z_dist}).
Therefore analogous to~(\ref{i1p_diff}) we get
\begin{equation}\label{j1_est}
|J_1-J^{(P)}_1| = \sum_{\frac{\eta_1 n}{2N} \leq k_1,k_2 \leq \frac{2\eta_2 n}{N}} |B_{l_1,l_2} - Poi(k_1;np_{l_1}) Poi(k_2;np_{l_2})| \Delta(k_1,q_{l_1}) \Delta(k_2,q_{l_2})
\end{equation}
where~\(q_{l_1}, q_{l_2}\) and~\(\Delta(.,.)\) are as defined in~(\ref{del_p_def}) and
\(B_{l_1,l_2} = B(k_1,k_2;n,p_{l_1},p_{l_2})\) is as defined in~(\ref{bin_dist2}).

Since~\(k_1\) and~\(k_2\) are both of the order of~\(\frac{n}{N},\) we get from~(\ref{b_poi_diff2}) that
\begin{equation}\nonumber
|B_{l_1,l_2} - Poi(k_1;np_{l_1})Poi(k_2;np_{l_2})| \leq \delta Poi(k_1,np_{l_1})Poi(k_2;np_{l_2}) \frac{n}{N^2}
\end{equation}
for some constant~\(\delta > 0\) not depending on~\(n,k_1,k_2,l_1\) or~\(l_2.\)
Using this in~(\ref{j1_est}) and arguing as in~(\ref{i1_temp_estt}) we then get
\begin{eqnarray}
|J_1-J^{(P)}_1| \leq C \mathbb{E}_0(T^{(P)}_{l_1}) \mathbb{E}_0(T^{(P)}_{l_2}) \left(\frac{n}{N^2}\right) \nonumber
\end{eqnarray}
for some constant~\(C >0.\) Using the upper bound~\(\mathbb{E}_0(T^{(P)}_{l_1}) \leq C_1 r_n \sqrt{\frac{n}{N}}\) for some constant~\(C_1\) not depending
on~\(l_1\) (see~(\ref{del_tn})), we then get
\begin{eqnarray}
|J_1-J^{(P)}_1| \leq C_2 \frac{r_n^2n^2}{N^3}. \label{j1_est2}
\end{eqnarray}
Substituting~(\ref{j1_est2}) and~(\ref{i12_est_fin2}) into~(\ref{t12_p1}) gives the final estimate in~(\ref{tl_diff2}). The middle
estimate in~(\ref{tl_diff2}) follows from the bounds for~\(\mathbb{E}T_l\) in~(\ref{del_tn_b}).~\(\qed\)

\emph{Proof of Lemma~\ref{cov_lemma}}:
Since the Poisson process~\({\cal P}\) is independent on disjoint subsets,
we have \[cov_0(T^{(P)}_{l_1},T^{(P)}_{l_2}) = \mathbb{E}_0(T^{(P)}_{l_1}T^{(P)}_{l_2}) - \mathbb{E}_0T^{(P)}_{l_1}\mathbb{E}_0T^{(P)}_{l_2} = 0.\]

Therefore write
\[|cov(T_{l_1},T_{l_2})| = |cov(T_{l_1},T_{l_2}) - cov_0(T^{(P)}_{l_1},T^{(P)}_{l_2})| \leq Z_1 + Z_2 + Z_3,\]
where
\[Z_1 = |\mathbb{E} T_{l_1} T_{l_2} - \mathbb{E}_0 T^{(P)}_{l_1} T^{(P)}_{l_2}| \leq C\left(\frac{r_n^2 n^2}{N^3}\right),\]
\[Z_2 = |\mathbb{E}^{(P)}_0T_{l_1}\mathbb{E}^{(P)}_0T_{l_2} - \mathbb{E}T_{l_1}\mathbb{E}T_{l_2}| \leq Z_3 + Z_4,\]
\[Z_3 = |\mathbb{E}_0T^{(P)}_{l_1} - \mathbb{E}T_{l_1}|\mathbb{E}_0T^{(P)}_{l_2} \leq C\left(\frac{r_n n^{3/2}}{N^{5/2}}\right)
\left(r_n\sqrt{\frac{n}{N}}\right) = C \left(\frac{r_n^2 n^2}{N^3}\right)\]
and similarly,
\[Z_4 = \mathbb{E}T_{l_1}|\mathbb{E}_0T^{(P)}_{l_2} - \mathbb{E}T_{l_2}| \leq  C\frac{r_n^2 n^2}{N^3},\]
for some constant~\(C > 0.\) The estimate for~\(Z_1\) follows from~(\ref{tl_diff2}) and the estimates for~\(Z_3\) and~\(Z_4\)
follow from~(\ref{tl_diff}) and the estimates for~\(\mathbb{E}T_l\) and~\(\mathbb{E}_0T^{(P)}_l\) in~(\ref{del_tn_b})
and~(\ref{del_tn}), respectively.~\(\qed\)


\setcounter{equation}{0}
\renewcommand\theequation{\thesection.\arabic{equation}}
\section{Proof of Theorem~\ref{tsp_thm}}\label{pf_tsp}
For~\(1 \leq  l \leq N,\) recall from~(\ref{tl_bd_sm}) that~\(T_l\) is the length of the minimum length cycle~\({\cal C}_l\) containing all the nodes
of~\(\{X_k\}\) contained in the square~\(S_l.\) Also we have from Section~\ref{intro} that~\(s_n\) denotes
the minimum distance between two squares in~\(\{S_l\}_{1 \leq l \leq N}.\) If the squares in~\(\{S_l\}\) are sufficiently far apart it is intuitive to expect that the overall minimum length cycle~\({\cal C}_{tot}\) containing all the nodes of~\(\{X_k\}\) is simply obtained by merging together the cycles~\({\cal C}_l.\) In other words, it is reasonable to expect that~\({\cal C}_{tot}\) ``covers" all nodes of a particular square before ``proceeding" to the next square. However, we give a small argument below to see that this is not necessarily true if the total number of nodes~\(n\) is large enough.

Suppose the intercity distance~\(s_n = 10r_n\) and~\(r_n = \sqrt{M\frac{\log{n}}{n}}\) for some large constant~\(M >0.\)
If all the~\(r_n \times r_n\) squares in Figure~\ref{sq_plc} are populated with nodes, then total number of squares~\(N\) satisfies~\[C_1\frac{n}{\log{n}} \leq \frac{1}{(20r_n)^2} \leq N \leq \left(\frac{1}{r_n}\right)^2 \leq C_2 \frac{n}{\log{n}}\]
for some constants~\(C_1,C_2 > 0.\) Condition~(\ref{N_est})
is therefore satisfied and so the estimates for the expected length of~\(T_l\) in Lemma~\ref{tl_lemma} hold.
From~(\ref{del_tn}) we therefore have that~\[\mathbb{E}T_l \geq C_3 r_n \sqrt{\frac{n}{N}} \geq C_4 r_n \sqrt{\log{n}}\]
for some constants~\(C_3, C_4 > 0.\) In other words, the expected total length of a cycle containing all the nodes of~\(S_l\)
is much larger than the intercity distance~\(s_n.\) Therefore it is quite possible that the cycle~\({\cal C}_{tot}\)
locally crosses between two squares \(s_n\) apart multiple times.

We now allow~\(s_n\) and~\(r_n\) to be general
as in the statement of the Theorem~\ref{tsp_thm} and show that the length~\(TSPC_n\) of the minimum length cycle~\({\cal C}_{tot}\)
is \emph{well approximated} by~\(\sum_{l=1}^{N} T_l.\)
\begin{Lemma}\label{tl_approx_lem}
The overall minimum length
\begin{equation}\label{tspc_up}
TSPC_n \leq \left(V_n + 2N(s_n+8r_n)\right) \ind(U_{tot}(n)) + 5\sqrt{n} \ind(U^c_{tot}(n)),
\end{equation}
where
\begin{equation}\label{vn_def1}
V_n := \sum_{l=1}^{N} T_l,
\end{equation}
\begin{equation}\label{vn_def}
U_{tot} = U_{tot}(n) := \bigcap_{l=1}^{N} U_l
\end{equation}
and~\(U_l\) is the event defined in~(\ref{ul_def}).
If the intercity distance~\(s_n > r_n \sqrt{2},\) then
\begin{equation}\label{tspc_low}
TSPC_n \geq V_n.
\end{equation}
\end{Lemma}

\emph{Proof of~(\ref{tspc_up})}: Suppose that the event~\(U_{tot}\) occurs and let~\({\cal C}_l\) be minimum length cycle
containing all the nodes in the square~\(S_l, 1 \leq l \leq N.\) Call the cycles~\(\{{\cal C}_l\}_{1 \leq l \leq N}\)
as \emph{small cycles}.
We construct a big cycle containing all the~\(n\) nodes by merging the small cycles~\({\cal C}_l\) together
iteratively, via a sequence of intermediate cycles~\(\{{\cal T}(i)\}_{1 \leq i \leq N}\) as follows.
Let~\({\cal T}(1) = {\cal C}_1\) so that the length of~\({\cal T}(1)\) is
\begin{equation}\label{lt1}
L({\cal T}(1)) = L({\cal C}_1) = T_1.
\end{equation}

To proceed with the iteration, recall from Section~\ref{intro} that the squares~\(\{S_l\}\) are well connected in the sense
that there exists a square in~\(\{S_j\}_{2 \leq j \leq N}\) at a distance~\(s_n\) from~\(S_1.\)
Without loss of generality, we assume that~\(S_i, 2 \leq i \leq N\) is
at a distance~\(s_n\) from some square~\(S_{q(i)} \in \{S_1,\ldots,S_{i-1}\}.\)

\begin{figure}[tbp]
\centering
\includegraphics[width=2in, trim= 80 370 100 190, clip=true]{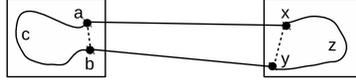}
\caption{Merging cycles~\({\cal T}(1) = {\cal C}_1 = acba\) and~\({\cal C}_2 = xzyx.\)}
\label{cyc_merg}
\end{figure}

Consider the small cycle~\({\cal C}_2\) containing all the nodes of~\(S_2.\) Remove any
edge~\(e_1\) from the intermediate cycle~\({\cal T}(1)\) and any edge~\(e_2\) from~\({\cal C}_2\)
and add ``cross edges"~\(f_1\) and~\(f_2\)
connecting the endvertices of~\(e_1\) and~\(e_2.\)
This is illustrated in Figure~\ref{cyc_merg} where the edges~\(e_1 = ab\)
and~\(e_2= xy\) are replaced by the edges~\(f_1 = ax\) and~\(f_2 = by.\)

The resulting intermediate cycle~\({\cal T}(2)\) satisfies the following properties with~\(i = 2\):\\
\((f1)\) The cycle~\({\cal T}(i)\) contains all the edges of the small cycles~\(\{{\cal C}_{j}\}_{1 \leq j \leq i}\)
not removed so far in the iteration process.\\
\((f2)\) The length
\begin{equation}\label{lt2}
L({\cal T}(i)) \leq L({\cal T}(i-1)) + 2(s_n + 8r_n) \leq \sum_{l=1}^{i}T_l + 2(i-1)(s_n + 8r_n).
\end{equation}
Property~\((f1)\) is true by construction and property~\((f2)\) is true
since the length of each added edge~\(f_i, i =1,2\) is no more than~\(s_n + 8r_n,\)
the sum of the distance between the squares~\(S_1\) and~\(S_2\) and the total perimeter of~\(S_1\) and~\(S_2.\)


Consider now a general iteration step~\(i \geq 3\)
where we need to merge the intermediate cycle~\({\cal T}(i)\) with the small cycle~\({\cal C}_{i+1}\)
containing all the nodes in the square~\(S_{i+1}.\) Recall that the square~\(S_{i+1}\)
is at a distance of~\(s_n\) from some square~\(S_{q(i)} \in \{S_1,\ldots,S_{i-1}\}.\)

Since the event~\(U_{tot}\) occurs, each square~\(S_l,1 \leq l \leq N\) contains at least
\[\frac{\eta_1 n}{2N} \geq \frac{\eta_1 M}{2} \log{n} \geq 8,\]
nodes of~\(\{X_k\}\) for all large~\(n,\) by~(\ref{n_N}). In particular,~\(S_{q(i)}\)
also contains at least~\(8\) nodes and so the small cycle~\({\cal C}_{q(i)}\)
contains at least~\(8\) edges.

The square~\(S_{i+1}\) is at a distance of~\(s_n\)
from~\(S_{q(i)}\) and so there are at most three squares in~\(\{S_j\}_{1 \leq j \leq i-1}\)
at a distance of~\(s_n\) from~\(S_{q(i)}.\) This means at most three edges
have been removed from the small cycle~\({\cal C}_{q(i)}\) in the iteration process so far
and so by property~\((f1),\) at least one edge~\(e_{q(i)}\) of~\({\cal C}_{q(i)}\) is still present in
the intermediate cycle~\({\cal T}(i).\)

Remove~\(e_{q(i)}\) and an edge from~\({\cal C}_{i+1}\) and add cross edges as before
to get the new cycle~\({\cal T}(i+1).\) Arguing as above, the new intermediate cycle~\({\cal T}(i+1)\)
also satisfies properties~\((f1)-(f2).\) Performing the above process for a total of~\(N-1\) iterations, we finally
obtain a big cycle~\({\cal C}_{fin}\) containing all the nodes~\(\{X_i\}_{1 \leq i \leq n},\)
whose length satisfies~
\begin{equation}\label{c_fin}
L({\cal C}_{fin}) \leq \sum_{i=1}^{N} T_l + 2(N-1)(s_n + 8r_n).
\end{equation}
Since the overall minimum length~\(TSPC_n \leq L({\cal C}_{fin})\) we obtain the upper bound~(\ref{tspc_up})
when~\(U_{tot}(n)\) occurs.

If the event~\(U_{tot}(n)\) does not occur, then we use the strips estimate (\ref{tsp_ab}) with~\(a = n\) and~\(b = 1\)
to get that the minimum length cycle~\(TSPC_n\) has a total length of at most~\(5\sqrt{n}.\)~\(\qed\)

\emph{Proof of~(\ref{tspc_low})}: For illustration we consider the case of two squares first.
Let~\({\cal Q}_1 = (v_1,\ldots, v_{k_1},v_1=:v_{k_1+1})\) be  the minimum length cycle containing all the nodes in~\(S_1\)
and let~\({\cal Q}_2 = (u_1,\ldots,u_{k_2},u_1 =: u_{k_2+1})\) be minimum length cycle containing all the nodes
in~\(S_2.\)
If~\({\cal C}_{tot}\) is the minimum length cycle containing all the nodes~\(\{v_j\} \cup \{u_j\},\) then
\begin{equation}\label{l_tot_est}
L\left({\cal C}_{tot}\right) \geq L\left({\cal Q}_1\right) + L\left({\cal Q}_2\right)
\end{equation}
where~\(L\left({\cal Q}_j\right), j =1,2\) is length of the cycle~\({\cal Q}_j\) as defined in~(\ref{len_cyc_def}).\\
\emph{Proof of~(\ref{l_tot_est})}: For a node~\(v \in \{v_j\} \cup \{u_j\},\) let~\(l\left(v,{\cal C}_{tot}\right)\)
be the sum of  length of the edges containing the node~\(v\) in the cycle~\({\cal C}_{tot}.\) Using~(\ref{len_cyc_def})
\begin{equation}\label{c_tot}
L\left({\cal C}_{tot}\right) = I_1 + I_2,
\end{equation}
where
\begin{equation}\label{i1_def}
I_1 = \frac{1}{2} \sum_{j =1}^{k_1} l\left(v_j,{\cal C}_{tot}\right) \text{ and } I_2 = \frac{1}{2} \sum_{j=1}^{k_2} l\left(u_j,{\cal C}_{tot}\right).
\end{equation}

\begin{figure}[tbp]
\centering
\includegraphics[width=2in, trim= 80 430 130 200, clip=true]{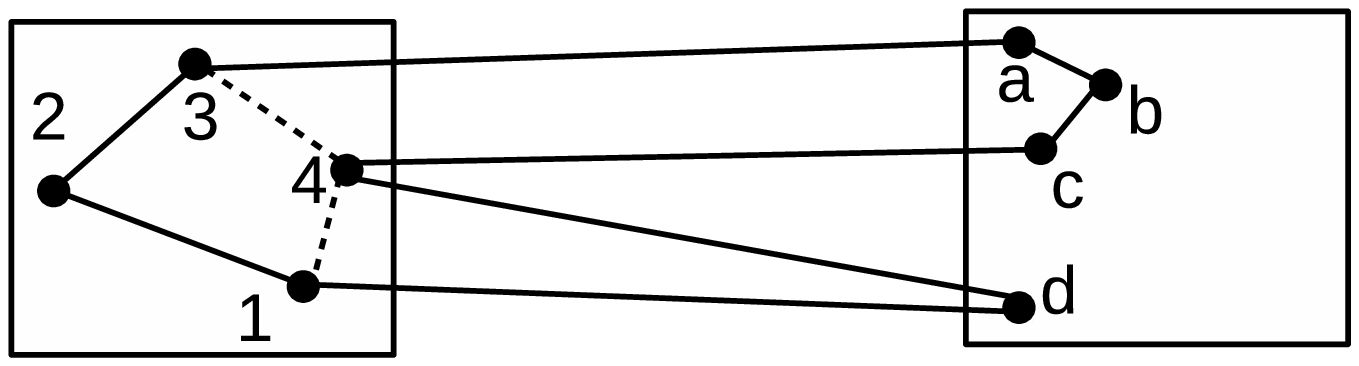}
\caption{Replace cycle~\({\cal C}_{tot} = 123abc4d1\) with the cycle~\({\cal C}_1 = 12341.\)}
\label{path_rep}
\end{figure}

To estimate~\(I_1,\) assume without loss of generality that the cycle~\({\cal C}_{tot}\) is of the form~
\begin{equation}\label{pv_rep}
{\cal C}_{tot} = (v_1,{\cal E}_1, v_2, {\cal E}_2,\ldots,{\cal E}_{k_1-1},v_{k_1}, {\cal E}_{k_1}, v_1 = v_{k_1+1}),
\end{equation}
where each~\({\cal E}_j\) is either empty or is a path containing only nodes of~\(\{u_j\}.\)
For~\(1 \leq j \leq k_1,\) replace the subpath~\({\cal E}_j\) of~\({\cal C}_{tot}\) with the edge~\((v_{j},v_{{j+1}}).\)
Let~\({\cal C}_1\) be the resulting cycle as shown in Figure~\ref{path_rep},
where~\(v_i\) is denoted by~\(i\) for~\(1 \leq i \leq 4.\)

For any fixed~\(1 \leq j \leq k_1\) the sum length of the edges containing~\(v_j\) as an endvertex is
less in the new cycle~\({\cal C}_1\) than in the original cycle~\({\cal C}_{tot}\) i.e.,
\begin{equation}\label{less_length}
l(v_j,{\cal C}_1) \leq l(v_j,{\cal C}_{tot})
\end{equation}
To see~(\ref{less_length}) is true, let~\(e_1\) and~\(e_2\) be the edges of~\({\cal C}_{tot}\)
containing~\(v_j\) as an endvertex in the original cycle~\({\cal C}_{tot}.\) Using the representation
of~\({\cal C}_{tot}\) in~(\ref{pv_rep}), we assume that the other endvertex of~\(e_1\) is either~\(v_{j-1}\)
or a node in~\(\{u_k\}.\) If~\(v_{j-1}\) is the other endvertex of~\(e_1,\) then~\(e_1\) is also
present in the new cycle~\({\cal C}_1.\) Else the length of~\(e_1\) is at least~\(s_n  > r_n \sqrt{2}\)
and~\(e_1\) is replaced by the edge~\(f_1 = (v_{j-1},v_j)\) in~\({\cal C}_1.\)
The length of~\(f_1\) is at most~\(r_n\sqrt{2}\) since both endvertices of~\(f_1\) lie within
the~\(r_n \times r_n\) square~\(S_1.\) A similar argument holds for the edge~\(e_2\) and so~(\ref{less_length}) is true.

Using~(\ref{less_length}) in~(\ref{i1_def}), we have
\begin{equation}
I_1  \geq \frac{1}{2} \sum_{j=1}^{k_1} l\left(v_j,{\cal C}_{1}\right) = L({\cal C}_1) \geq L\left({\cal Q}_1\right),\label{ip_def2}
\end{equation}
since~\({\cal Q}_1\) is the minimum length cycle containing all the nodes~\(\{v_j\}.\) An analogous argument obtains that~\(I_2 \geq L\left({\cal Q}_2\right)\) and so from~(\ref{c_tot}), we get~(\ref{l_tot_est}). The argument for the general case is analogous.~\(\qed\)

We use Lemma~\ref{tl_approx_lem} to prove Theorem~\ref{tsp_thm}.
From Lemma~\ref{tl_approx_lem},
we have that the overall minimum length~\(TSPC_n\) is bounded above and below by the sum of
the local minimum lengths~\(\sum_{l=1}^{N}T_l.\)
From the bounds on~\(\mathbb{E}T_l\) in~(\ref{del_tn}) of Lemma~\ref{tl_lemma},
we have that~\(\sum_{l=1}^{N} \mathbb{E}T_l\) is of the order of~\(N r_n \sqrt{\frac{n}{N}}  = r_n \sqrt{nN} = b_n\)
as defined in~(\ref{bn_def}). We therefore study the convergence of~\(\frac{TSPC_n}{b_n}.\)
We henceforth fix~\(M > 0\) large so that~(\ref{cov_tl_est}) of Lemma~\ref{cov_lemma} holds.\\\\
\emph{Proof of~(\ref{conv_tsp_prob}) in Theorem~\ref{tsp_thm}}:
From the upper and lower bounds~(\ref{tspc_up}) and~(\ref{tspc_low}) in Lemma~\ref{tl_approx_lem},
we have that
\begin{eqnarray}
\frac{1}{b_n}(V_n - \mathbb{E}V_n)  - \Delta_n \leq \frac{1}{b_n}\left(TSP_n - \mathbb{E} TSP_n\right) \leq \frac{1}{b_n}(V_n - \mathbb{E}V_n) + \Delta_n \label{key_est2}
\end{eqnarray}
where~\(V_n = \sum_{l=1}^{N} T_l\) is as defined~(\ref{vn_def1}) and
\begin{equation}\label{delta_n_def}
\Delta_n = \frac{2(N-1)(s_n+8r_n)}{b_n}\ind(U_{tot}(n)) + \frac{5\sqrt{n}}{b_n}\ind(U^c_{tot}(n)) \nonumber.
\end{equation}

The variance of~\(V_n\) satisfies
\begin{equation}\label{var_vn}
var(V_n) \leq C \frac{r_n^2 n^2}{N} = C b_n^2 \left(\frac{n}{N^2}\right)
\end{equation}
for some constant~\(C > 0\) and all~\(n\) large and  since~\(\frac{n}{N^2} \longrightarrow 0\) (see~(\ref{N_est})),
we get that
\begin{equation}\label{vn_prob_conv}
\frac{1}{b_n}\left(V_n - \mathbb{E}V_n\right) \longrightarrow 0 \text{ in probability }
\end{equation}
as~\(n \rightarrow \infty.\) Also
\begin{equation}\label{del_n_conv}
\Delta_n \longrightarrow 0 \text{ a.s.}
\end{equation}
as~\(n \rightarrow \infty.\) This proves~(\ref{conv_tsp_prob}) and we prove~(\ref{var_vn}) and~(\ref{del_n_conv}) separately below.

\emph{Proof of~(\ref{var_vn})}: Write
\begin{eqnarray}
var(V_n) &=& \sum_{l} var(T_l) + \sum_{l_1,l_2} cov(T_{l_1},T_{l_2}) \nonumber\\
&\leq& \sum_{l} \mathbb{E}T_l^2  + \sum_{l_1,l_2} cov(T_{l_1},T_{l_2}), \label{var_vn_est1}
\end{eqnarray}
where~\(cov(X,Y) = \mathbb{E}XY - \mathbb{E}X\mathbb{E}Y.\)  Using~(\ref{del_tn}) of Lemma~\ref{tl_lemma} to
estimate~\(\mathbb{E}T_l^2\) we get~
\begin{equation}\label{f_t}
\sum_{l=1}^{N} \mathbb{E} T_l^2 \leq N C_1 \left(r_n \sqrt{\frac{n}{N}}\right)^2 = C_1r_n^2 n
\end{equation}
for some constant~\(C_1 > 0.\) Similarly using estimate~(\ref{cov_tl_est}) of Lemma~\ref{cov_lemma} for the covariance,
we get
\begin{equation}\label{s_t}
\sum_{l_1,l_2} cov(T_{l_1},T_{l_2})\leq N^2 \left(C_2 \frac{r_n^2 n^2}{N^3}\right) = C_2 \frac{r_n^2 n^2}{N}.
\end{equation}
for some constants~\(C > 0.\) Substituting~(\ref{f_t}) and~(\ref{s_t}) into~(\ref{var_vn_est1}), we get
\[var(V_n) \leq C_1 r_n^2 n + C_2 \frac{r_n^2 n^2}{N} = \frac{r_n^2 n^2}{N}\left(C_1 \frac{N}{n} + C_2\right).\]
Since~\(\frac{N}{n} \leq \frac{1}{M\log{n}} \leq 1\) for all~\(n\) large~(see~(\ref{n_N})), we get that~\(var(V_n) \leq C_3\frac{r_n^2 n^2}{N}\)
for some positive constant~\(C_3\) and for all~\(n\) large.

\emph{Proof of~(\ref{del_n_conv})}: From~(\ref{delta_n_def}) and the fact that~\(r_n < r_n\sqrt{2} < s_n\) (see statement of the Theorem),
we get~
\begin{equation}\label{delta_n}
0 \leq \Delta_n \leq \frac{18N s_n}{b_n}  + \frac{5\sqrt{n}}{b_n} \ind(U^c_{tot}(n))
\end{equation}
and so
\begin{equation}\label{delta_n2}
0 \leq \limsup_n \Delta_n \leq \limsup_n \frac{5\sqrt{n}}{b_n} \ind(U^c_{tot}(n)),
\end{equation}
since \(\frac{Ns_n}{b_n} \longrightarrow 0\) as~\(n \rightarrow \infty\) by the statement of the Theorem.
From the estimate for the event~\(U_l\) in~(\ref{ul_def}),
\begin{equation}\label{u_tot_est}
\mathbb{P}(U^c_{tot}(n)) \leq \sum_{l=1}^{N} \mathbb{P}(U^c_l) \leq N \exp\left(-C\frac{n}{N}\right),
\end{equation}
for some constant~\(C > 0.\) Using the fact that~\(\frac{n}{N} \geq M\log{n}\) (see~(\ref{n_N})), we get
\begin{equation}\label{u_tot_est2}
\mathbb{P}(U^c_{tot}(n)) \leq \frac{n}{M\log{n}} \frac{1}{n^{MC}} \leq \frac{1}{n^2},
\end{equation}
provided~\(M > 0\) is large. Fixing such an~\(M,\) we have from Borell-Cantelli lemma that~\(\mathbb{P}(\limsup_n U^{c}_{tot}(n))=0\) and so a.s. \(\ind(U_{tot}^c(n)) = 0\) for all large~\(n.\) From~(\ref{delta_n2}), we therefore get~(\ref{del_n_conv}).~\(\qed\)

\emph{Proof of~(\ref{exp_tspc}) in Theorem~\ref{tsp_thm}}: Recalling that~\(V_n = \sum_{i=1}^{N} T_l\)
from~(\ref{vn_def1}), we use Lemma~\ref{tl_approx_lem}
to get
\begin{equation}\label{vn_est_tl}
\mathbb{E}V_n \leq \mathbb{E} TSPC_n \leq \mathbb{E}V_n + b_n\mathbb{E}\Delta_n,
\end{equation}
where~\(\Delta_n\) satisfies (see~(\ref{delta_n}))
\begin{equation}\label{edel_n}
\mathbb{E}\Delta_n \leq \frac{18Ns_n}{b_n}  + \frac{5\sqrt{n}}{b_n} \mathbb{P}(U_{tot}^c(n)) \leq 18 + \frac{5\sqrt{n}}{b_n} \mathbb{P}(U_{tot}^c(n)),
\end{equation}
since~\(\frac{Ns_n}{b_n} \longrightarrow 0\) as~\(n \rightarrow \infty\) (see statement of the Theorem).
Using~(\ref{u_tot_est2}) for estimating the probability of the event~\(U_{tot}\) we get
\begin{equation}\label{u_tot_bn}
\sqrt{n}\mathbb{P}(U^{c}_{tot}(n)) \leq \frac{\sqrt{n}}{n^2} \leq \sqrt{\frac{M\log{n}}{n}} \leq r_n
\end{equation}
for all~\(n\) large, where the final inequality is true by the condition for~\(r_n\) in~(\ref{N_est}). 

On the other hand~\(b_n = r_n \sqrt{nN} \geq r_n\) and so we get from~(\ref{u_tot_bn}) that
\begin{equation}\label{u_tot_bn2}
\frac{5\sqrt{n}}{b_n}\mathbb{P}(U^{c}_{tot}(n))  \leq 5
\end{equation}
and using~(\ref{u_tot_bn2}) in~(\ref{edel_n}) we get~\(\mathbb{E}\Delta_n \leq 23\) and so from~(\ref{vn_est_tl}),
\begin{equation}\label{vn_est_tl2}
\mathbb{E}V_n \leq \mathbb{E} TSPC_n \leq \mathbb{E}V_n + 23b_n.
\end{equation}

To estimate~\(\mathbb{E}V_n\) use the bounds for~\(\mathbb{E}T_l\) in~(\ref{del_tn}) of Lemma~\ref{tl_lemma} to get
\begin{equation}\label{vn_estb4}
C_1 b_n = N\left(C_1 r_n \sqrt{\frac{n}{N}}\right) \leq \mathbb{E}V_n \leq N\left(C_2 r_n \sqrt{\frac{n}{N}}\right) = C_2 b_n
\end{equation}
for some constants~\(C_1,C_2  >0.\) From~(\ref{vn_estb4}) and~(\ref{vn_est_tl2}), we get the bounds for~\(\mathbb{E}TSPC_n\) in~(\ref{exp_tspc}).~\(\qed\)

\emph{Proof of~(\ref{eq_tsp2}) of Theorem~\ref{tsp_thm}}:
We consider Poissonization and recall the Poisson process~\({\cal P}\) on the squares~\(\{S_l\}_{1 \leq l \leq N},\)
defined on the probability space~\((\Omega_0, {\cal F}_0, \mathbb{P}_0)\) (see paragraph prior to~(\ref{poi_dist})).
Analogous to~\(TSPC_n\) defined~(\ref{min_weight_cycle}), let~\(TSPC^{(P)}_n\) denote the length of the minimum length cycle containing all
the nodes of the Poisson process~\({\cal P}.\) Recall from~(\ref{tl_def_p}) that~\(T^{(P)}_l\) denotes the length
of the minimum length cycle containing all the nodes of~\({\cal P}\) in the square~\(S_l.\)

Analogous to~(\ref{tspc_low}), we have that
if the intercity distance~\(s_n > r_n \sqrt{2},\) then
\begin{equation}\label{tspc_low_p}
TSPC^{(P)}_n \geq V^{(P)}_n = \sum_{l=1}^{N} T_l^{(P)}.
\end{equation}

Define the event~\[E^{(P)}_l = \left\{T^{(P)}_l \geq \delta_4 r_n \sqrt{\frac{n}{N}}\right\},\] where~\(\delta_4\) is the constant
in~(\ref{tn_prob_est}) of Lemma~\ref{tl_lemma_poiss}. Since the Poisson process
is independent on disjoint sets, the events~\(E^{(P)}_l\) are independent
and each occurs with probability at least~\(\delta_5,\) by~(\ref{tn_prob_est}).
If
\begin{equation}\label{fp_sum}
F^{(P)}_{sum} := \sum_{l=1}^{N} \ind(E^{(P)}_l)
\end{equation}
then~\(\mathbb{E}_0 \left(F^{(P)}_{sum}\right) \geq \delta_5 N\)
and from the standard Chernoff bound estimate for sums of independent Bernoulli
random variables (see Corollary \(A.1.14,\) pp. 312 of Alon and Spencer~(2008))
we also have
\begin{equation}\label{v_est}
\mathbb{P}_0\left(F^{(P)}_{sum} \geq C_1 N\right)
\geq 1- e^{-2C_2 N}
\end{equation}
for some positive constants~\(C_1\) and~\(C_2.\) If~\(F^{(P)}_{sum} \geq C_1N,\) then by~(\ref{fp_sum}),
the sum
\[\sum_{l=1}^{N} T^{(P)}_l \geq C_1 N \left(\delta_4 r_n \sqrt{\frac{n}{N}}\right) = C_3 b_n\]
for some constant~\(C_3 > 0\) and so from~(\ref{tspc_low_p}),
\begin{equation}\label{poi_t_up}
\mathbb{P}_0(TSPC^{(P)}_n \geq C_3 b_n) \geq 1- e^{-2C_2 N}
\end{equation}
for all~\(n\) large.

To convert the probability estimates to the Binomial process, let~\[A_P = \{TSPC^{(P)}_n \geq C_3 b_n\}, A = \{TSPC_n \geq C_3 b_n\}\]
and use the dePoissonization formula
\begin{equation}\label{de_poiss_ax}
\mathbb{P}(A) \geq 1- D \sqrt{n} \mathbb{P}(A^c_P)
\end{equation}
for some constant~\(D > 0\) and~(\ref{poi_t_up}) to get that
\begin{equation}\label{poi_t_up2}
\mathbb{P}(TSPC_n \geq C_3 b_n) \geq 1- D\sqrt{n} e^{-2C_2 N} = 1-e^{-  \alpha_N},
\end{equation}
where
\[\alpha_N = 2C_2N - \log{D} - \frac{1}{2}\log{n} \geq C_2N\] for all~\(n\) large, since~\(N \geq \sqrt{n}\) for all~\(n\) large~(see~(\ref{n_N})).
This proves~(\ref{eq_tsp2}) and it only remains to prove~(\ref{de_poiss_ax}).

To prove~(\ref{de_poiss_ax}), let~\(N_P\) denote the random number of nodes of~\({\cal P}\) in all the squares~\(\cup_{j=1}^{N}S_j\) so
that~\(\mathbb{E}_0 N_P = n\) and~\(\mathbb{P}_0(N_P=n) = e^{-n}\frac{n^{n}}{n!} \geq \frac{D_1}{\sqrt{n}}\) for
some constant~\(D_1 > 0,\) using the Stirling formula. Given~\(N_P = n,\) the nodes of~\({\cal P}\)
are i.i.d.\ with distribution~\(g_N\) as defined in~(\ref{gn_def}); i.e.,
\[\mathbb{P}_0(A_P^c|N_P = n)  = \mathbb{P}(A^c)\] and so
\[\mathbb{P}_0(A_P^c) \geq \mathbb{P}_0(A_P^c|N_P = n) \mathbb{P}_0(N_P = n) =
\mathbb{P}(A^c) \mathbb{P}_0(N_P = n) \geq \mathbb{P}(A^c)\frac{D_1}{\sqrt{n}},\] proving~(\ref{de_poiss_ax}).~\(\qed\)

\emph{Proof of~(\ref{eq_tsp3}) of Theorem~\ref{tsp_thm}}:
As in the proof of~(\ref{eq_tsp2}) above, we consider the Poisson process~\({\cal P}\)
on the squares~\(\{S_l\}_{1 \leq l \leq N}\) defined in the paragraph prior to~(\ref{poi_dist}).
As before, let~\(TSPC^{(P)}_n\) denote the length of the minimum length cycle containing all
the nodes of the Poisson process~\({\cal P}.\) Recall from~(\ref{tl_def_p}) that~\(T^{(P)}_l\) denotes the length
of the minimum length cycle containing all the nodes of~\({\cal P}\) in the square~\(S_l.\)

Analogous to~(\ref{tspc_up}), we have
\begin{equation}\label{tspc_up_p}
TSPC^{(P)}_n \leq \left(V^{(P)}_n + 2(N-1)(s_n+8r_n)\right) \ind(U^{(P)}_{tot}(n)) + 5\sqrt{n} \ind(U^{(P)}_{tot}(n))^c,
\end{equation}
where
\begin{equation}\label{vn_def1_p}
V^{(P)}_n := \sum_{l=1}^{N} T^{(P)}_l,
\end{equation}
\begin{equation}\label{vn_def_p}
U^{(P)}_{tot} = U^{(P)}_{tot}(n) := \bigcap_{l=1}^{N} U^{(P)}_l
\end{equation}
and~\(U^{(P)}_l = \{\frac{\eta_1 n}{2N} \leq N^{(P)}_l \leq \frac{2\eta_2 n}{N}\}\) is the event defined in~(\ref{ul_def_p}).
Recall that~\(N^{(P)}_l\) is the total number of nodes of~\({\cal P}\) inside the square~\(S_l.\)

Suppose now that the event~\(U_{tot}^{(P)}(n)\) occurs so that
\begin{equation}\label{tspc_ph}
TSPC_n^{(P)} \leq V_n^{(P)}  + 2(N-1)(s_n + 8r_n) = \sum_{l=1}^{N} T^{(P)}_l + 2(N-1)(s_n+8r_n).
\end{equation}
Since~\(U^{(P)}_l \supseteq U_{tot}^{(P)}\) occurs for every~\(1 \leq l \leq N,\)
we use the strips estimate~(\ref{tsp_ab}) with~\(a = \frac{2\eta_2 n}{N}\) and~\(b = r_n\)
to get that the corresponding minimum length~\(T_l^{(P)} \leq 5b\sqrt{a} \leq C r_n \sqrt{\frac{n}{N}}\)
for some constant~\(C > 0\) and for every~\(1 \leq l \leq N.\)
Thus \[V^{(P)}_n = \left(\sum_{l=1}^{N} T^{(P)}_l\right) \leq C b_n\]
and from~(\ref{tspc_ph})
we therefore get
\begin{equation}\label{tspc_ph1}
TSPC_n^{(P)} \leq C b_n + 2(N-1)(s_n + 8r_n) \leq Cb_n + 18Ns_n \leq (C+1)b_n,
\end{equation}
for all~\(n\) large. The second inequality in~(\ref{tspc_ph1}) is true since~\(r_n < r_n \sqrt{2} < s_n.\)
The final inequality in~(\ref{tspc_ph1}) is true since~\(\frac{Ns_n}{b_n} \longrightarrow 0\) and
so~\(\frac{Ns_n}{b_n} \leq \frac{1}{18}\) for all~\(n\) large.

Summarizing, we have that if the event~\(U_{tot}^{(P)}\) occurs,
then the overall minimum length~\(TSPC_n^{(P)} \leq C_1 b_n\)
for some constant~\(C_1 > 0.\)
To evaluate~\(\mathbb{P}(U^{(P)}_{tot}),\) use the estimate~(\ref{ul_est_p}) for the event~\(U^{(P)}_l\) to get
\begin{equation}\label{u_tot_p_est}
\mathbb{P}_0(U^{(P)}_{tot}) \geq 1- N\exp\left(-2C\frac{n}{N}\right)
\end{equation}
for some constant~\(C > 0.\) Thus
\begin{equation}\label{tutu}
\mathbb{P}_0\left(TSPC_n^{(P)} \leq C_1 b_n \right) \geq \mathbb{P}(U^{(P)}_{tot}) \geq 1-N\exp\left(-2C\frac{n}{N}\right).
\end{equation}

To convert the probabilities to the Binomial process, we again use
the dePoissonization formula~(\ref{de_poiss_ax})
to get that
\begin{equation}\label{tutu_bin}
\mathbb{P}\left(TSPC_n \leq C_1 b_n \right) \geq 1-DN\sqrt{n}\exp\left(-2C\frac{n}{N}\right) = 1-e^{-\delta_N},
\end{equation}
where~\(D > 0\) is as in~(\ref{de_poiss_ax}) and
\begin{equation}\label{del_yt}
\delta_N = 2C \frac{n}{N} - \log{D}-\log{N} - \frac{1}{2}\log{n}.
\end{equation}
Since~\(\frac{n}{N} \geq M\log{n}\) for all~\(n\) large (see~(\ref{n_N})), we get
\[\log{D} + \log{N} + \frac{1}{2}\log{n} \leq \log{D} + \log\left(\frac{n}{M\log{n}}\right) + \frac{1}{2}\log{n} \leq 2\log{n} \leq C\frac{n}{N},\]
provided~\(M >0\) is large. Fixing such an~\(M\) we get that~\(\delta_{N} \geq C\frac{n}{N}\) and so~(\ref{eq_tsp3}) follows
from~(\ref{tutu_bin}).~\(\qed\)


\setcounter{equation}{0}
\renewcommand\theequation{\thesection.\arabic{equation}}
\section{Proof of Theorem~\ref{var_tsp_thm}}\label{pf_tsp_var}
We need preliminary estimates regarding the
change in length of the minimum length cycle upon adding or deleting a single node.

Let~\(X_1,\ldots,X_{n+1}\) be~\(n+1\) random nodes distributed according to the density~\(f\)
in the unit square~\(S.\) For~\(1 \leq j \leq n+1,\) let~\({\cal D}_j\)
denote the minimum length cycle containing all the nodes~\(\{X_k\}_{1 \leq k \neq j \leq n+1}\) with length
\begin{equation}\label{ldj}
L({\cal D}_j) = TSP(X_1,\ldots,X_{j-1},X_{j+1},\ldots,X_{n+1};S),
\end{equation}
where~\(TSP(.;.)\) is as defined in~(\ref{min_cyc}).
For future use, we estimate lengths of edges in~\({\cal D}_j.\)

Divide the unit square~\(S\) into~\(2Aw_n \times 2Aw_n\)
squares~\(\{W^{(1)}_i\}_{1 \leq i \leq N_W},\) each of side length~\(2Aw_n\) where
\begin{equation}\label{a_def}
\frac{1}{n^{1/6}} \leq w_n := \frac{1+c_n}{n^{1/6}} \leq \frac{2}{n^{1/6}}, \;\;\;A := \left(\frac{3}{\epsilon_1}\right)^{\frac{1}{3}}
\end{equation}
and~\(\epsilon_1  >0\) is as in~(\ref{f_eq}). The term~\(c_n \in (0,1)\) is chosen such that~\(\frac{1}{2Aw_n}\) is an integer
for all~\(n\) large. For~\(1 \leq i \leq N_W,\) let~\(W^{(2)}_i\) be the bigger square
with same centre as~\(W^{(1)}_i\) but with side length~\(4Aw_n.\)
For~\(1 \leq j \leq n+1\) and~\(1 \leq i \leq N_W,\) let~\(F_{j}(i)\) be the event
there exists an edge~\(e_j(i) \in {\cal D}_j\) with both endvertices in the bigger square~\(W^{(2)}_i\)
and let
\begin{equation}\label{f_tot_def}
F_{tot}(n+1) := \bigcap_{j=1}^{n+1} \bigcap_{i=1}^{N_W} F_j(i).
\end{equation}

The following Lemma is used in the proof of Theorem~\ref{var_tsp_thm}.
\begin{Lemma}\label{fji_lemm}
We have that
\begin{equation}\label{f_tot_est}
\mathbb{P}\left(F_{tot}(n+1)\right) \geq 1-\exp\left(-Cn^{2/3}\right)
\end{equation}
for some constant~\(C>0\) and for all~\(n \geq 3.\)
\end{Lemma}
\emph{Proof of Lemma~\ref{fji_lemm}}:
We first perform some preliminary computations.
Fix~\(1 \leq j \leq n+1\) and~\(1 \leq i \leq N_W.\) Using~(\ref{f_eq}) and
the fact that~\(w_n \geq n^{-\frac{1}{6}}\) (see~(\ref{a_def})),
the average number of
nodes of~\(\{X_k\}_{1 \leq k \neq j \leq n+1}\) in the square~\(W^{(1)}_i\) is
\[ n\int_{W^{(1)}_i} f(x) dx \geq n \epsilon_1 (2Aw_n)^{2}  \geq 4A^2 \epsilon_1 n^{2/3},\]
where~\(\epsilon_1 > 0\) is as in~(\ref{f_eq}).
Let~\(Z_j(i)\) denote the event that the square~\(W^{(1)}_i\)
contains at least~\(2A^2 \epsilon_1 n^{2/3}\) nodes of~\(\{X_k\}_{1 \leq k \neq j \leq n+1}.\)
By standard Binomial estimates
(see Corollary~\(A.1.14,\) pp. 312 of Alon and Spencer (2008)) and the fact that~\(w_n^2 \geq n^{-\frac{1}{3}}\)
(see~(\ref{a_def})), we get
\begin{equation}\label{vj_est}
\mathbb{P}\left(Z_j(i)\right) \geq 1-\exp(-C_1 nw_n^2) \geq 1-\exp\left(-C_2n^{2/3}\right)
\end{equation}
for some positive constants~\(C_1\) and~\(C_2.\)

If
\begin{equation}\label{z_tot_def}
Z_{tot}(n+1) := \bigcap_{j=1}^{n+1} \bigcap_{i=1}^{N_W} Z_j(i),
\end{equation}
then we have from~(\ref{vj_est}) that
\begin{equation}\label{z_tot_est2}
\mathbb{P}\left(Z_{tot}(n+1)\right) \geq 1-(n+1)N_W \exp\left(-C_2 n^{2/3}\right).
\end{equation}
The total number of squares is
\begin{equation}\label{nw_def}
N_W = \left(\frac{1}{2Aw_n}\right)^2 \geq D n^{1/3}
\end{equation}
for some constant~\(D > 0\) using~\(w_n \geq n^{-\frac{1}{6}}\) (see~(\ref{a_def})) and so
we get from~(\ref{z_tot_est2}) that
\begin{equation}\label{z_tot_est}
\mathbb{P}\left(Z_{tot}(n+1)\right) \geq 1-\exp\left(-C_3 n^{2/3}\right)
\end{equation}
for some constant~\(C_3 > 0.\)

The estimate~(\ref{z_tot_est}) and the following property imply Lemma~\ref{fji_lemm}.\\
\((f1)\) If the event~\(Z_{tot}(n+1)\) occurs, then for every~\(1 \leq j \leq n+1\) and~\(1 \leq i \leq N_W,\)
there exists an edge~\(e_j(i) \in {\cal D}_j\) with both endvertices in the bigger square~\(W^{(2)}_i.\)\\
\emph{Proof of~\((f1)\)}: Suppose~\(Z_{tot}(n+1)\) occurs and suppose that
the node~\(X_j\) is present in the square~\(W^{(1)}_i.\)
Let~\(\{Y_k\}_{1 \leq k \leq q} \subset \{X_k\}_{1 \leq k \neq j \leq n+1}\)
be the other nodes present in the square~\(W^{(1)}_i.\)
Since the event~\(Z_j(i) \supseteq Z_{tot}(n+1)\) occurs,
\begin{equation}\label{q_est}
q \geq 2A^2 \epsilon_1 n^{2/3}.
\end{equation}


For~\(1 \leq k \leq q,\) let~\(e_k(1)\) and~\(e_k(2)\) be the edges containing the node~\(Y_k\) as an endvertex
in the cycle~\({\cal D}_j.\) If no edge of~\({\cal D}_j\) has both its endvertices inside the bigger square~\(W^{(2)}_i,\)
then all the edges~\(\{e_k(1),e_k(2)\}_{1 \leq k \leq q}\) are distinct and each such edge
has length at least~\(Aw_n,\) since it must cross the annulus\\\(W^{(2)}_i \setminus W^{(1)}_i.\)
Therefore if~\(l(Y_k,{\cal D}_j)\) is the sum of length of the edges containing~\(Y_k\)
as an endvertex in the cycle~\({\cal D}_j,\) then~\(l(Y_k, {\cal D}_j) \geq 2Aw_n.\)

From~(\ref{len_cyc_def}) we therefore have that the
total length of~\({\cal D}_j\) is
\begin{equation}\label{len_cyc_dj}
L({\cal D}_j) \geq \frac{1}{2} \sum_{k=1}^{q} l(Y_k,{\cal D}_j) \geq q\cdot Aw_n \geq 2A^3\epsilon_1 n^{2/3} w_n.
\end{equation}
Using the fact that~\(w_n \geq n^{-\frac{1}{6}}\) (see~(\ref{a_def})) we then get that
\begin{equation}\label{ldj_cont}
L({\cal D}_j) \geq 2A^3 \epsilon_1 \sqrt{n} \geq 6\sqrt{n},
\end{equation}
by our choice of~\(A\) in~(\ref{a_def}).

But using the strips estimate~(\ref{tsp_ab}) with~\(a = n\) and~\(b = 1,\)
we have that the length of the cycle~\({\cal D}_j\) is at most
\begin{equation}\nonumber
L({\cal D}_j) \leq 5 \sqrt{n}
\end{equation}
and this contradicts~(\ref{ldj_cont}).~\(\qed\)

The above Lemma allows us to estimate the variance of the
length of the minimum length cycle.\\
\emph{Proof of~\ref{var_tsp_est_main} of Theorem~\ref{var_tsp_thm}}:
We use the martingale difference method and
for~\(1 \leq j \leq n+1,\) let~\[{\cal F}_j = \sigma\left(X_1,\ldots,X_j\right)\] denote the sigma field
generated by the random variables~\(X_1,\ldots,X_j.\) Defining the martingale difference
\begin{equation}\label{di_diff}
D_j = \mathbb{E}(TSP_{n+1} | {\cal F}_j) - \mathbb{E}(TSP_{n+1} | {\cal F}_{j-1}),
\end{equation}
we have~\[TSP_{n+1} -\mathbb{E}TSP_{n+1} = \sum_{j=1}^{n+1} D_j\] and so by the martingale property
\begin{equation}\label{var_tsp_di}
var(TSP_{n+1})  = \mathbb{E}\left(\sum_{j=1}^{n+1} D_j\right)^2 = \sum_{j=1}^{n+1} \mathbb{E}D_j^2.
\end{equation}
There is a constant~\(C> 0\) such that
\begin{equation}\label{di_sec}
\max_{1 \leq j  \leq n+1}\mathbb{E}D_j^2 \leq \frac{C}{n^{1/3}}
\end{equation}
for all~\(n \geq 1\) and this proves~(\ref{var_tsp_est_main}).

\emph{Proof of~(\ref{di_sec})}: We first rewrite~\(D_j\) in a more convenient form.
Let~\(\omega = (x_1,\ldots,x_{n+1})\) and~\(\omega' = (y_1,\ldots,y_{n+1})\)
be two vectors in~\((\mathbb{R}^2)^{n+1}.\)
Defining~\[\omega_j = (x_1,\ldots,x_j,y_{j+1},\ldots,y_{n+1})\] for~\(1 \leq j \leq n+1\)
and using Fubini's theorem, we get
\begin{equation}\label{di_est2}
|D_j| = \left|\int (T(\omega_j) - T(\omega_{j-1}))f(y_j)\ldots f(y_n)dy_{j}\ldots dy_{n+1} \right|  \leq H_j \nonumber\\
\end{equation}
where
\begin{equation}
H_j := \int |T(\omega_j) - T(\omega_{j-1})|f(y_j)\ldots f(y_n)dy_{j}\ldots dy_{n+1} \label{di_est}
\end{equation}
and~\(T(\omega_t), t = j,j-1\) is the length of the minimum length cycle containing all the nodes in~\(\omega_t.\)\\

Let~\(F_{tot}(n+1)\) be the event defined in~(\ref{f_tot_def})
and write
\begin{equation}\label{d_decom}
H_j  = I_1  + I_2,
\end{equation}
where
\begin{eqnarray}
I_1 &=& \int |T(\omega_j) - T(\omega_{j-1})|\ind(\omega_j \in F_{tot}(n+1))\ind(\omega_{j-1} \in F_{tot}(n+1)) \nonumber\\
 &&\;\;\;\;\;\;\;\;f(y_j)\ldots f(y_{n+1})dy_{j}\ldots dy_{n+1} \label{i1_def_f}
\end{eqnarray}
and~\(I_2 = I_1 - H_j.\)

There is a positive constant~\(C  >0\) such that
\begin{equation}\label{i1_tot_est}
\mathbb{E}I^2_1 \leq \frac{C}{n^{1/3}}  \;\;\;\;\;\text{ and } \;\;\;\;\; \mathbb{E}I_2^2 \leq \exp\left(-C n^{2/3}\right)
\end{equation}
and so using~\(H_j^2 =(I_1+I_2)^2 \leq 2(I_1^2 + I_2^2)\) we get
\[\mathbb{E}H_j^2 \leq 2\left(\frac{C}{n^{1/3}}  + \exp\left(-C n^{2/3}\right)\right) \leq  \frac{3C}{n^{1/3}}\]
for all~\(n\) large. This proves~(\ref{di_sec}).

We obtain the estimates for~\(I_1\) and~\(I_2\) in~(\ref{i1_tot_est}) separately below.\\
\underline{Estimate for~\(I_1\)}: Let~\({\cal D}_j\) be the minimum length cycle containing
all the nodes~\(\{x_k\}_{1 \leq k \leq j-1} \cup \{y_k\}_{j+1 \leq k \leq n+1}.\)
If~\(L({\cal D}_j)\) is the length of~\({\cal D}_j,\) then for~\(t \in \{j-1,j\}\)
\begin{eqnarray}
|T(\omega_t) - L({\cal D}_t)| \ind(\omega_t \in F_{tot}(n+1)) \leq 4Aw_n \sqrt{2} \label{t_om_td}
\end{eqnarray}
and so from~(\ref{t_om_td}),~(\ref{i1_def_f}) and triangle inequality, we have
\begin{equation}\label{i1_estr}
I_1 \leq 8Aw_n\sqrt{2} \text{   and so   } \mathbb{E}(I^2_1) \leq C_1 w_n^2 \leq \frac{C_2}{n^{1/3}}
\end{equation}
for some constants~\(C_1,C_2\) since~\(w_n \leq \frac{2}{n^{1/6}}\) (see~(\ref{a_def})).\\
\emph{Proof of~(\ref{t_om_td})}: We prove for~\(t=j\) and an analogous analysis holds
for~\(t = j-1.\) By monotonicity~(\ref{tsp_ab}), we have that~\(T(\omega_j) \geq L({\cal D}_j).\)
Also, since~\(\omega_j \in F_{tot}(n+1),\) every square~\(W^{(2)}_k, 1 \leq k \leq N_W\)
of side length~\(4Aw_n\) defined prior to Lemma~\ref{fji_lemm} contains an edge of~\({\cal D}_j.\)
Suppose the ``new" node~\(x_j\) belongs to the square~\(W^{(1)}_i.\) Since there is an edge~\(e \in {\cal D}_j\)
having both its endnodes~\(z_1,z_2\) inside~\(W^{(1)}_i,\) we remove~\(e\) and add the edges~\((z_1,x_j)\)
and~\((x_j,z_2)\) to form a new cycle containing all the nodes of~\(\omega_j.\) The total length of the two edges added
is at most~\(4Aw_n \sqrt{2}.\) This implies that~\(T(\omega_j) \leq L({\cal D}_j) +4Aw_n \sqrt{2},\) proving~(\ref{t_om_td}).

\underline{Estimate for~\(I_2\)}: Every edge within the unit square~\(S\)
has length at most~\(\sqrt{2}\) and any cycle containing all the~\(n+1\) nodes of~\(\omega_t\)
has~\(n+1\) edges. Therefore~\(T(\omega_t) \leq (n+1)\sqrt{2}\) for~\(t \in \{j-1,j\}.\) Thus from the definition
of~\(I_2\) in~(\ref{d_decom}), we get~\(I_2 \leq J_1 + J_2,\) where
\[J_1  = (n+1)\sqrt{2} \int \ind(\omega_j \notin F_{tot}(n+1))f(y_j)\ldots f(y_n)dy_j\ldots dy_n\]
and
\[J_2 = (n+1)\sqrt{2} \int \ind(\omega_{j-1} \notin F_{tot}(n+1))f(y_j)\ldots f(y_n)dy_j\ldots dy_n.\]
Using Cauchy-Schwarz inequality,
\[J_1^2 = 2(n+1)^2 \left(\mathbb{E}(\ind(F^c_{tot}(n+1))|{\cal F}_j)\right)^2 \leq 2(n+1)^2 \mathbb{E}(\ind(F^c_{tot}(n+1))|{\cal F}_j)\] and
similarly~\[J_2^2 \leq 2(n+1)^2 \mathbb{E}\left(\ind(F^c_{tot}(n+1)) | {\cal F}_{j-1}\right).\]

Since~\(I_2^2  = (J_1 +J_2)^2 \leq 2(J_1^2 + J_2^2)\) and
\(\mathbb{E}(\mathbb{E}(X|{\cal F}_j)|{\cal F}_{j-1}) = \mathbb{E}(X | {\cal F}_{j-1}),\)
we get that
\begin{equation}\nonumber
\mathbb{E}(I_2^2|{\cal F}_{j-1}) \leq 4(n+1)^2\mathbb{P}(F^c_{tot}(n+1) | {\cal F}_{j-1})
\end{equation}
and therefore that
\begin{equation}\label{i2_estr}
\mathbb{E}I_2^2 \leq 4(n+1)^2 \mathbb{P}\left(F^c_{tot}(n+1)\right) \leq 4(n+1)^2\exp\left(-2C n^{2/3}\right) \leq \exp\left(-Cn^{2/3}\right)
\end{equation}
for some constant~\(C > 0\) and for all~\(n\) large, using~(\ref{f_tot_est}).~\(\qed\)

\emph{Proof of~(\ref{exp_tsp_u}) and~(\ref{eq_tsp1_u}) of Theorem~\ref{var_tsp_thm}}:
The upper bound for~\(\mathbb{E} TSP_n\) in~(\ref{exp_tsp_u})
is obtained from the strips estimate~(\ref{tsp_ab}) with~\(a =n\) and~\(b = 1.\)
This also proves~(\ref{eq_tsp1_u}).

For the lower bound in~(\ref{exp_tsp_u}), we argue as follows.  For~\(1 \leq i \leq n,\) let~\(d(X_i,\{X_j\}_{1 \leq j \neq i \leq n})\) denote the minimum distance of node~\(X_i\) from all other nodes~\(\{X_j\}_{1 \leq j \neq i \leq n}.\) The TSP length~\(TSP_n\) then satisfies
then~\(TSP_n \geq \sum_{i=1}^{n} d(X_i,\{X_j\}_{1 \leq j \neq i \leq n})\) and so
\begin{equation}\label{tsp_er}
\mathbb{E}TSP_n \geq n\mathbb{E}d(X_1,\{X_j\}_{2 \leq j \leq n}).
\end{equation}
Analogous to the proof of~(\ref{x2_est}), we have that~\[\mathbb{E}d(X_1,\{X_j\}_{2 \leq j \leq n}) \geq \frac{C}{\sqrt{n}}\]
for some constant~\(C > 0\) not depending on the choice of~\(i\) and so from~(\ref{tsp_er}) we get~(\ref{eq_tsp2_u}).~\(\qed\)

\emph{Proof of~(\ref{eq_tsp2_u})}: Divide the unit square~\(S\) into~\(r_n \times r_n\) squares~\(\{S_l\}_{1 \leq l \leq N}\) placed~\(s_n\)
apart as in Figure~\ref{sq_plc} with~\(r_n\) and~\(s_n\) as follows:
\begin{equation}\label{rn_sn_def2}
r_n^2 = \frac{M\log{n} + c_n}{n} \text{ and } s_n^2 = \frac{2M\log{n} + d_n}{n}
\end{equation}
where~\(c_n \in (0,1)\) and~\(d_n \in (4,5)\) are such that~\(\frac{1-r_n}{r_n+s_n}\)
is an integer. With this choice of~\(r_n\) and~\(s_n,\) the number of~\(r_n \times r_n\) squares~\(N\)
and the scaling factor~\(b_n\) defined in~(\ref{bn_def}) satisfy
\begin{equation}\label{n_def_tt}
C_1 nr_n^2 \leq C_2 \frac{n}{\log{n}} \leq  N \leq C_3\frac{n}{\log{n}} \leq C_4 nr_n^2
\end{equation}
and
\begin{equation}\label{bn_Est}
C_5 \sqrt{n} \leq b_n = r_n \sqrt{nN} \leq C_6 \sqrt{n}
\end{equation}
for some positive constants~\(\{C_i\}_{1 \leq i \leq 6}.\)

Let~\({\cal C}_n\) denote the minimum length cycle
containing all the nodes of~\(\{X_k\}_{1 \leq k \leq n}\) present
in all the~\(r_n \times r_n\) squares~\(\{S_j\}_{1 \leq j \leq N}.\)
If~\(L({\cal C}_n)\) denotes the length of~\({\cal C}_n,\) then by monotonicity~(\ref{tsp_ab})
we have that
\begin{equation}\label{tsp_ineq1}
TSP_n \geq L({\cal C}_n)
\end{equation}
and  since the term~\(s_n > r_n \sqrt{2}\) strictly (see~(\ref{rn_sn_def2})),
we have from~(\ref{tspc_low}) that
\begin{equation}\label{tsp_ineq2}
L({\cal C}_n) \geq \sum_{l=1}^{N} T_l,
\end{equation}
where~\(T_l, 1 \leq l \leq N\) is the minimum length cycle containing all the nodes of~\(\{X_k\}\)
in the square~\(S_l.\)

Estimates for~\(\mathbb{E}T_l\) in Lemma~\ref{tl_lemma}
and estimates for~\(\mathbb{E}_0T^{(P)}_l,\)
the Poissonized process, in Lemma~\ref{tl_lemma_poiss} hold in this case as well.
Moreover if~\(M > 0\) is large in~(\ref{rn_sn_def2}), then the covariance
estimate in Lemma~\ref{cov_lemma} holds as well.
For illustration, we prove the lower bound for~\(\mathbb{E}T_l\) here.
From~(\ref{f_eq}), any node of~\(\{X_k\}_{1 \leq k \leq n}\) is present in the square~\(S_l\)
with probability~
\[ D_1\frac{n}{N} \leq \epsilon_1 r_n^2 \leq q_l = \int_{S_l} f(x)dx \leq \epsilon_2 r_n^2 \leq D_2 \frac{n}{N}\]
for some positive constants~\(D_1\) and~\(D_2,\) using~(\ref{n_def_tt}). The estimates
for~\(q_l\) are  analogous to the estimates for~\(p_l\) in~(\ref{pl_def}). Arguing as in the proof of~(\ref{del_tn_b})
we then get that~\(\mathbb{E}T_l \geq C r_n \sqrt{\frac{n}{N}}.\)

Arguing as in the proof of~(\ref{eq_tsp2}), we get
\[\mathbb{P}\left(TSP_n \geq L({\cal C}_n) \geq C_4 b_n\right) \geq 1-e^{-C_5 N}\]
for some positive constants~\(C_4,C_5.\) Finally, using~(\ref{n_def_tt}) and~(\ref{bn_Est}) to estimate~\(b_n\) and~\(N\)
we get~(\ref{eq_tsp2_u}).~\(\qed\)





\subsection*{Acknowledgement}
I thank Professors Rahul Roy and Federico Camia for crucial comments and for my fellowships.

\bibliographystyle{plain}

\begin{thebibliography}{10}







\bibitem{alon} N. Alon and J. Spencer. (2008).
\newblock{\em The probabilistic method}.
\newblock{Wiley}.

\bibitem{arora} S. Arora. (1998).
\newblock{ PTAS for Euclidean Traveling Salesman and
Other Geometric Problems}.
\newblock{\em Journal of the ACM}, \textbf{45}, pp. 753--782.

\bibitem{beard} J. Beardwood, J. H. Halton and J. M. Hammersley. (1959).
\newblock{The shortest path through many points}.
\newblock{\em Proceedings Cambridge Philosophical Society}, \textbf{55}, pp. 299--327.

\bibitem{cook} W. Cook. (2011).
\newblock{\em In pursuit of the Traveling Salesman: Mathematics at the Limits of Computation}.
\newblock{Princeton University Press}.

\bibitem{gutin} G. Gutin and A. P. Punnen. (2006).
\newblock{\em The traveling salesman problem and its variations}.
\newblock{Springer}.

\bibitem{karp} M. Karpinski, M. Lampis and R. Schmied (2015).
\newblock{New Inapproximability bounds for TSP}.
\newblock{\em Journal of Computer and System Sciences}, \textbf{81}, 1665–-1677.

\bibitem{pint} C-M. Pintea, P. C. Pop and C. Chira. (2017)
\newblock{The generalized traveling salesman problem solved with
ant algorithms}.
\newblock{\em Complex Adaptive Systems Modeling}, \textbf{5}, available
at \emph{https://doi.org/10.1186/s402940017-0048-9}.

\bibitem{snyd} L. V. Snyder and M. S. Daskin. (2006).
\newblock{A random-key genetic algorithm for the generalized traveling
salesman problem}.
\newblock{\em European Journal of Operational Research}, \textbf{174}, 38--53.

\bibitem{steele2} J. M. Steele. (1981).
\newblock{Subadditive Euclidean functionals and nonlinear growth in geometric probability}.
\newblock{\em Annals of Probability}, \textbf{9}, pp. 365--376.

\bibitem{steele} J. M. Steele. (1993).
\newblock{Probability and Problems in Euclidean Combinatorial Optimization}.
\newblock{\em Statistical Science}, \textbf{8}, pp. 48--56.

\bibitem{stin} S. Steinerberger. (2015).
\newblock{New bounds for the Traveling Salesman constant}.
\newblock{\em Advances in Applied Probability}, \textbf{47}, pp. 27--36.

\bibitem{vazi} V. V. Vazirani. (2001).
\newblock{\em Approximation algorithms}.
\newblock{Springer Verlag}.








\end{thebibliography}

\end{document}